\documentclass[review]{elsarticle}

\date{}
\usepackage{lineno,hyperref,epstopdf}
\usepackage{amsmath,amsthm,amssymb}
\usepackage{caption}
\usepackage{pdflscape}
\usepackage{times}
\usepackage{natbib}
\usepackage{tabularx}
\usepackage{array}
\usepackage{xcolor}
\usepackage{bigstrut}
\usepackage{booktabs}
\usepackage{graphicx}
\usepackage{float}
\usepackage{multirow}
\usepackage{url}

\modulolinenumbers[5]

\newtheorem{theorem}{Theorem}

\newtheorem{corollary}{Corollary}
\newtheorem{proposition}{Proposition}
\newtheorem{remark}{Remark}

\begin{document}

\begin{frontmatter}

\title{On IPW-based estimation of conditional average treatment effect}
\author[BNU]{Niwen Zhou}

\author[BNU,HKBU]{Lixing Zhu\corref{mycorrespondingauthor}\corref{thank} }
\cortext[thank]{
			The authors gratefully acknowledge two grants from the University Grants Council of Hong Kong  and a NSFC grant (NSFC11671042). \hspace{.2cm}}
\cortext[mycorrespondingauthor]{Corresponding author}
\ead{lzhu@hkbnu.edu.hk}

\address[BNU]{School of Statistics, Beijing Normal University, Beijing, China}

\address[HKBU]{Department of Mathematics, Hong Kong Baptist University, Hong Kong, China}

\begin{abstract}
 The research in this paper gives a systematic investigation on the asymptotic behaviours of four inverse probability weighting (IPW)-based estimators for conditional average treatment effect, with nonparametrically, semiparametrically, parametrically estimated and true propensity score, respectively.
To this end,  we first pay a particular attention to semiparametric dimension reduction structure such that  we can well study the semiparametric-based estimator that can well alleviate the curse of dimensionality and greatly avoid model misspecification. We also derive some further properties of  existing  estimator with nonparametrically estimated propensity score. According to their asymptotic variance functions, the studies reveal the general ranking of their asymptotic efficiencies; in which scenarios the asymptotic equivalence can hold; the critical roles of the affiliation of the given covariates in the set of arguments of propensity score, the bandwidth and kernel selections. The results show
an essential difference from the  IPW-based (unconditional) average treatment effect(ATE).
The numerical studies indicate that for high-dimensional paradigms, { the semiparametric-based estimator performs
well in general {whereas nonparametric-based estimator}, even sometimes, parametric-based estimator, is more affected by dimensionality.}  Some numerical studies are carried out to examine their performances. A real data example is
analysed for illustration.
\end{abstract}


\end{frontmatter}

\section{Introduction}

Treatment effects have been widely analyzed by economists and statisticians in diverse fields. In this paper,
we focus on estimating treatment effect under the potential outcomes framework {and the unconfoundedness assumption} with binary treatment.
Let $D=0,1$ mean that the individual does not receive or receives
treatment and the response $Y$ be the corresponding potential outcome as $Y(0)$ or
$Y(1)$. {To conveniently identify the quantities measuring treatment effects, the unconfoundedness assumption in \citep{Rosenbaum1983} is generally considered, that is, the assignment to  treatment is independent of the potential outcomes given a $k$-dimensional vector $X$ of covariates, i.e.
\begin{eqnarray}{\label{uncon_ass}}
(Y(0),Y(1))\perp D\mid X.
\end{eqnarray}}
Further, we in this paper consider
the dimension of $X$ to be fixed throughout this paper, but in some cases it can be high.\footnote{ 
Although the word "high dimension" is usually conjunct with $k$ being divergent with sample size in recent years, when we say $X$ is of high dimension in this paper, it only means $X$ contains many but fixed number of covariates.
For ease of explanation, we still use the word "high dimension" whenever no confusion will be caused.
}
As $Y(0)$ and $Y(1)$ cannot be simultaneously observed for any individual, the observed outcome can be written as $Y=DY(1)+(1-D)Y(0).$ Since estimating the $i$-th individual treatment
effect $(Y_i(1)-Y_i(0))$ is unrealistic, an important trend in the literature turns to
estimate the average treatment effect ($ATE$):~  $\mu=E(Y(1)-Y(0)).$
See for instance {\cite{Rosenbaum1983} and \cite{HIR2003}}.

{Recently, there is an increasing interest in estimating conditional (or heterogeneous) average treatment effects: $CATE(X)=E(Y(1)-Y(0)\mid X),$ which is designed to reflect how treatment effects vary across different subpopulations. Note that even thought receiving a treatment may have no effect on outcomes for the overall population, i.e. $ATE=0$, the treatment can still be effective for a subpopulation defined by specific observable characteristics, i.e. for some $x$ such that $CATE(x)\neq 0.$ Thus heterogeneous treatment effects are more informative and can  play  important roles in  personalized medicine or policy intervention. Most of existing  estimation methods for the heterogeneous treatment effects are conditional on the full set of variables, $X$, see e.g. \cite{crump2008}, \cite{wager2018}, where the multivariate variable $X$ are designed to make the unconfoundedness assumption plausible.
After 2015, researchers consider to estimate more general conditional/heterogeneous treatment effects, in which the conditioning covariates $Z$ with $Z$ being a subset of covariates, i.e.
$$X=(Z^{\top}, U^{\top})^{\top}\in
R^{l}\times R^{m},k=l+m<\infty.$$
 See e.g. \cite{Abrevaya2015} and \cite{Lee2017}. Note that treatment effects conditioning on a subset of $X$, rather than the high dimensional covariates $X$, can provide desirable flexibility and can help making policy decision.

 Based on the assumption (\ref{uncon_ass}), \cite{Abrevaya2015} used the inverse probability weighting ($IPW$)-based method, which is  popularly used in literature \citep{robins1994}, to estimate $$CATE(Z)=E[Y(1)-Y(0)\mid Z]$$ when the
propensity score function is estimated parametrically ($IPW$-$P$) and
nonparametrically ($IPW$-$N$).  \cite{Abrevaya2015} gave a deep
investigation on the asymptotic properties of the estimators.
 There are two main conclusions in \cite{Abrevaya2015}: one is $IPW$-$N$  can be asymptotically
more efficient than $IPW$-$P$ in the sense that the asymptotic variance function of $IPW$-$N$ can be uniformly smaller than that of $IPW$-$P$, the another is the asymptotic variance function of $IPW$-$P$ is equal to that of $IPW$-$O$ which is defined as the
oracle estimator with the true propensity score.
{It is noteworthy that the last conclusion is different from that of IPW-type ATE estimators, because the IPW-type ATE estimator based on parametrically estimated propensity score can be more efficient than the one with true propensity score.}

As is known,  to make the unconfoundedness assumption be plausible, it is often the case that we need to include
many covariates  in the analysis. {Thus we say $X \in R^k$ is of high dimension with $k<\infty$.}
In this case, on one hand, it is often not easy to choose a
parametric specification that can sufficiently capture all the important
nonlinear and interaction effects to have $IPW$-$P$. On the other hand,  any nonparametric estimation of
propensity score clearly suffers from the curse of dimensionality and then
$IPW$-$N$ does not work any more.

Therefore, in this paper,  we suggest a semiparametric IPW-based $CATE(Z)$ estimation procedure to simultaneously alleviate  the propensity score misspecification problem and particularly
the curse of
dimensionality. To this end, we consider a semiparametric dimension reduction structure of the propensity score and the unconfoundedness assumption~\eqref{uncon_ass} can have a dimension reduction version. It is worth pointing out that the general nonparametric structure can be regarded as a special case of the dimension reduction structure we consider with an orthonormal projection matrix of full rank.   
We will call the estimator $IPW$-$S$ and give the details about the model setting and the estimation procedure in the next section.

For theoretical development, we will give the asymptotically linear
representation and asymptotic normality of $IPW$-$S$. We will also give
 some
further properties of existing $IPW$-$N$ in \cite{Abrevaya2015}. Based on the theoretical studies, we give a systematic comparison
 on the asymptotic
efficiency amongst $IPW$-$O$, $IPW$-$P$, $IPW$-$S$ and
 $IPW$-$N$.

Combining the results of \cite{Abrevaya2015}
and the further properties of $IPW$-$N$ we derive in this paper,
the comparison reveals some very interesting and important phenomena.
Specifically, letting $A\preceq B$ mean that the asymptotic variance of estimator $A$  is not greater than that of estimator $B$ and $A\cong B$ stand for that $A$ has the same asymptotic variance function as  $B$, we have the following observations in theory.

  First,  in general
 $IPW\mbox{-}N \preceq IPW\mbox{-}S \preceq IPW\mbox{-}P \cong IPW\mbox{-}O$.

  Second,  the affiliation of $Z$ to the set of arguments of the propensity score plays an important role in the
     asymptotic efficiency of $IPW$-$S$ and $IPW$-$N$. That is,
     when $Z$ is a subset of arguments of the propensity score,   $
     IPW\mbox{-}S  \preceq IPW\mbox{-}P  \cong IPW\mbox{-}O$ and  $
     IPW\mbox{-}N  \preceq IPW\mbox{-}P  \cong IPW\mbox{-}O$, otherwise,
 $IPW\mbox{-}N\cong IPW\mbox{-}S
\cong IPW\mbox{-}P \cong IPW\mbox{-}O$.
 Note that this newly found phenomenon provides a deep insight into the
performances of $IPW$-$S$ and $IPW$-$N$, which is also useful in practice.

  Third, when the propensity score function is smooth enough,  then
    even in general cases we can also have the asymptotic equivalence by carefully choosing the bandwidths and using high order kernel
    functions: $IPW\mbox{-}N \cong IPW\mbox{-}S \cong IPW\mbox{-}P \cong IPW\mbox{-}O$. This also gives
    us a better understanding for the asymptotic performance of different
    estimators. Of course, this part mainly serves as a theoretical
    exploration. For practical use, we would have no interest to wilfully choose
    those kernel function and bandwidths, which are very difficult to
    implement and make the estimator with worse performance. But it reminds
    the researchers that a ``good" estimator of the propensity score would not be helpful for the performance of
    the $CATE$ estimator.

  Fourth, owing to the dimension reduction structure of $p(X)$, the
    requirements for bandwidths and the order of kernel function for
    $IPW$-$S$ are much milder than those for $IPW$-$N$. Thus when the dimension
    is high, even though $IPW$-$N$ has the superior efficiency in theory,
    $IPW$-$S$ is preferable.

The rest of the paper is organized as follows. In Section 2, we first introduce the
estimation procedure for $IPW$-$S$. Also we investigate its
asymptotic properties  and the theoretical comparisons between the four $CATE$
estimators. Section 3 contains some numerical studies to examine the
performance of the $CATE$ estimators. In Section 4, we apply the CATE estimators to analyse a
real data set for illustration.
Section~5 contains some conclusions and a further discussion.
The regularity conditions are listed in Appendix and all the technical proofs are relegated to Supplementary Materials to save space.

\section{Semiparametric estimation procedure and asymptotic properties}

\subsection{ Preliminary of estimation}
Assume that covariates $X=(Z^{\top},U^{\top})^{\top}$ are absolutely
continuous, under the unconfoundedness assumption (\ref{uncon_ass}), recall that $CATE$ function
$\tau(z)$ can be rewritten as \begin{eqnarray}\label{def_cate}
\tau(z)=E\left[\frac{DY}{p(X)} -\frac{(1-D)Y}{1-p(X)} \mid Z=z\right], Z\in R^l.
\end{eqnarray}

If $p(X)$ is given, we can estimate $\tau(z)$ immediately via the
Nadaraya-Watson kernel method by regarding $ \frac{DY}{p(X)} -\frac{(1-D)Y}{1-p(X)}$ as response:
\begin{eqnarray*}\label{orcal_e}
\hat \tau_O(z)= \bigg({ \sum_{i=1}^n
\left[ \frac{D_iY_i}{ p(X_i)} -\frac{(1-D_i)Y_i}{1-
 p( X_i)} \right]K_{h}(Z_i-z)}\bigg)\Big /
{ \sum_{i=1}^nK_{h}(Z_i-z)}.
\end{eqnarray*}
Here $K(\cdot)$  is a multivariate kernel
 function, $K_{h}(u)= {h^{-l}}
K\left( {u}/{h}\right)$ and $l=\mbox{dim}(Z)$. This $CATE$ estimator is
$IPW$-$O$ we mentioned before.

{Based on existing results for nonparametric estimation, it is easy to derive the asymptotic distribution of $IPW$-$O$ which will be used as the benchmark to make comparisons among all estimators studied in this paper.
\begin{proposition} Suppose the conditions (C1)-(C4) in Appendix are satisfied, the following
statements  hold for each point $z$ in the support of $Z$:
\begin{eqnarray*}
\sqrt{nh^l}\left(\hat \tau_O(z)-\tau(z)\right)
\stackrel{\mathcal{D}}{\longrightarrow}N\left(0, \frac{\parallel
K \parallel_2^2\sigma_{O}^2(z)}{f(z)}
\right).
\end{eqnarray*}
Here $\sigma_{O}^2(z)=E\left(\left[ \frac{ DY}{ p(X)}- \frac{(1-D)Y}{1-
 p(X)}-\tau(z) \right]^2\mid Z=z\right).$
\end{proposition}
}

When $p(X)$ is an unknown function, we then
first estimate  $p(X)$ to define a final $CATE$ estimator $\hat \tau(z).$  We propose the estimator under  semiparametric structure below.

\subsection { Semiparametric estimation for conditional average treatment effect}

{Assume the propensity score has a semiparametric dimension reduction structure:
\begin{eqnarray}{\label{con_1}}
p(X)=q(V^{\top}X),
\end{eqnarray}
where both the function $q(\cdot)$} and the $r$ projection directions in $V$ are unknown with $V$ being a $k\times r$ orthonormal matrix. {It is noteworthy that  this structure is general, which covers the structures of some important semiparametric models such as single-index models.}
From the definition of
propensity score, (\ref{con_1}) implies  the indicator $D$ depends on $X$ through the
projected variable $V^{\top}X$. Thus, we can use the following conditional
independence to present the above semiparametric structure:
\begin{eqnarray}\label{conind}
D\perp X \mid V^{\top} X.
\end{eqnarray}

It follows that $(Y(0),Y(1))\perp D\mid V^{\top}X.$  We call the intersection of all $V$'s satisfying the above independence the central subspace,
see \cite{Li1991}. Usually $V$ can only be identified up to a rotation matrix $C$. That is, $V^*=V \times C$ can be identified.  As this identification issue does not affect the related estimation of $p(X)$, we then still use $V$ without confusion. Relevant references are \cite{Luo2017} and \cite{Ma2019}.
 This is a dimension reduction framework, so that the corresponding estimation could be less
affected
by the curse of dimensionality. For such a dimension reduction structure, we can also consider variable selection as \cite{Ma2019} did. But as this is not a focus of this paper, we then just work on this model and assume the existence of consistent estimation later on.

If we  postulate that  the information about $D$ from $X$ can be completely
captured by $r$ linear combinations $V^{\top}X$ of $X$ with $r\ll k$, the propensity score can be estimated by replaced the original  $X$ with $V^{\top}X$ . That is, we can use lower
dimensional kernel function $\mathcal{H}(u)$ to get a nonparametric
estimator $\hat q(\hat V^{\top }X)$ of $q(V^{\top}X)
=E(D\mid
V^{\top}X)$,
\begin{eqnarray}\label{semi_q}
\hat q(\hat V^{\top} X_i)=\frac{  \sum_{j\neq i}^nD_i
\mathcal{H}_{h_2}({\hat V^{\top}X_j-\hat V^{\top}X_i})} {
\sum_{j\neq i}^n\mathcal{H}_{h_2}({\hat V^{\top}X_j-\hat V^{\top}X_i})}.
\end{eqnarray}
where  $h_2$ is the bandwidth,
$\mathcal{H}_{h_2}(u)= {h_2^{-r}} \mathcal{H}\left({u}/{{h_2}}\right)$
and $\hat V$ is a consistent estimator derived by a sufficient dimension reduction
method. {There are several methods available in the literature,
such as inverse regression methods in
 \citep{Cook2002}
and minimum average variance
estimation(MAVE) in \citep{xia2002,xia2007} . }

 Recall that
 $CATE$ can be  rewritten as (\ref{def_cate}).
Thus, based on $\hat q(\hat V^{\top} X_i)$, the $IPW$-$S$ of  $\tau(z)$  is defined as
\begin{eqnarray}\label{est_scate}
\hat \tau_S(z)=  { \sum_{i=1}^n
\left[ \frac{D_iY_i}{\hat q(\hat V^{\top} X_i)}  - \frac{(1-D_i)Y_i}{1-\hat
 q(\hat V^{\top} X_i)} \right]K_{h}(Z_i-z)} \Big /
{ \sum_{i=1}^nK_{h}(Z_i-z)}.
\end{eqnarray}

  Since both $Z$ and $V^{\top}X$ are low-dimension random vectors, $\hat \tau_S(z)$ can well alleviate the propensity score misspecification problem and
the curse of
dimensionality simultaneously.

In the next section, we investigate the
 asymptotic properties of $\hat \tau_S(z)$ and derive some further properties of existing $IPW$-$N$ under certain regularity conditions.

\subsection{Asymptotic properties for IPW-S}
  Denote
 $|A|$ as the cardinality of set $A$. We first give some notations.
\begin{itemize}
\item[$\bullet$] $\mathcal{W}=(X,D,Y)$ and the observation data $\mathcal{W}_i=(X_i,D_i,Y_i)_{i=1}^n$
    are the independent copies of $\mathcal{W}$;

\item[$\bullet$] $m_{j}(V^{\top}X) ={E[Y(j)\mid V^{\top}X]},j=0,1,$ and
    $K_i=K( {(Z_{i} -z)}/{h})$;

\item[$\bullet$] $\psi(p(V^{\top}X),\mathcal{W}) = {DY}/\{q(V^{\top}X)\}
    - {(1-D)Y}/\{ 1-q(V^{\top}X) \}  $;

\item[$\bullet$] $\psi^*(q(V^{\top}X),\mathcal{W}) =  [D\{Y-m_1(V^{\top}
    X)\} ]/\{ q(V^{\top} X)\}- [(1-D)\{Y-m_0(V^{\top} X)\} ]/\{1- q(V^{\top} X)\}
    +m_1(V^{\top} X)-m_0(V^{\top} X) $.

 \item[$\bullet$]
  For two vectors $A$ and $B$, we use intersection notation $A\cap B$ to write, without confusion, as all components that are contained in both $A$ and $B$.  $|A\cap B|=t$  stands for the number of components in the intersection of $A$ and $B$.  Particularly, when $t=0$,
    $A\cap
     B=\varnothing$, and $t=|A|$ implies
     $A\cap B=A$.

\end{itemize}
Both $\psi(q(V^{\top}X),\mathcal{W})$ and $\psi^*(q(V^{\top}X),\mathcal{W})$
are the central parts of influence function for $IPW$-$S$.

\begin{theorem}{\label{sem_them}}
Suppose all the conditions in Appendix are satisfied, the following
statements  hold for each point $z$ in the support of $Z$:
\begin{itemize}
 \item[(1)]When $|Z \cap V^{\top}X|=t<l$ with
    $s_2[2-l/(l-t)]+l>0$, the asymptotically linear
    representation is
 \begin{align*}
\sqrt{nh^l} \left(\hat \tau_S(z)-\tau(z)\right)=\frac{1}{\sqrt{nh^l}f(z)} \sum_{i=1}^{n}[\psi(p(V^{\top}X_i),\mathcal{W}_i)-\tau(z)]K_i+o_p(1)
\end{align*}
 and the asymptotic distribution of $\hat \tau_S(z)$  is
\begin{align*}
\sqrt{nh^l} \left(\hat \tau_S(z)-\tau(z)\right)
\stackrel{\mathcal{D}}{\longrightarrow}N\left(0,\Sigma_S(z)\right).
\end{align*}
\item[(2)]When $|Z \cap V^{\top}X|=l$,   the asymptotically linear
    representation is
 \begin{align*}
\sqrt{nh^l} \left(\hat \tau_S(z)-\tau(z)\right)=\frac{1}{\sqrt{nh^l}f(z)} \sum_{i=1}^{n}[\psi^*(q(V^{\top}X_i),\mathcal{W}_i)
-\tau(z)]K_i+o_p(1)
\end{align*}
 and the asymptotic distribution of $\hat \tau_S(z)$ is
 \begin{eqnarray*}
\sqrt{nh^l} \left(\hat \tau_S(z)-\tau(z)\right)
\stackrel{\mathcal{D}}{\longrightarrow}N\left(0, \Sigma_S^*(z)
\right).
\end{eqnarray*}
\end{itemize}

Here $s_2$ is the order of $\mathcal{H}(\cdot)$, $\Sigma_S(z)= {\parallel
K \parallel_2^2\sigma_{S}^2(z)}\Big /{f(z)}
, \Sigma_S^*(z)= {\parallel
K \parallel_2^2\sigma_{S}^{*2}(z)}\Big /{f(z)}$ with
$\sigma_{S}^{2}(z) =E[\{\psi(q(V^{\top}X),\mathcal{W})-\tau(z)\}^2\mid
Z=z]$ and $\sigma_{S }^{*2}(z)=E[\{\psi^*(q(V^{\top}X),\mathcal{W})-\tau(z)\}^2\mid
Z=z]$.
\end{theorem}

\begin{remark}\label{ncate_noz}
 These results  show a very interesting and somewhat unexpected phenomenon that the asymptotic behaviors
of $\hat \tau_S(z) $ also depend on whether some of elements of $Z$ belong to  $V^{\top}X.$
 Recall that $|Z \cap V^{\top}X|=t$ means
 $t$ elements of $Z$ also are $t$ linear combinations of
$V^{\top}X$, i.e. we can rewrite $V^{\top}X=(Z_1,\cdots,Z_{t},(\tilde
V^{\top}X)^{\top})^{\top}$ with $V^{\top}=\left(\begin{array}{cc}
\tilde I_{(t)\times k} \\
\tilde V^{\top}
\end{array}\right)_{r\times k}$ and
$\tilde I_{(t)\times k}=\left(\begin{array}{c|c}I_{t\times t}
&0\end{array}\right).$
Here $I_{t\times t}$ is an identity
matrix, $\tilde V^{\top}$ is a $(r-t)\times k$ matrix. The asymptotic behaviours with $t=l$ and any  $0\leq t< l$ are very different. A natural question is whether we can, if possible, choose a dimension reduced vector  $V^{\top}X$ such that $IPW$-$S$ works best. The question is related to $IPW$-$P$ and $IPW$-$N$, we will have some detailed discission in Subsection~3.3 below.

\end{remark}

{Next, we present the estimators for $\Sigma_S(z)$ and $\Sigma^*_S(z)$  under $|Z\cap V^{\top}X|<l$ and $|Z\cap V^{\top}X|=l$ respectively as
 \begin{align}
 \hat \Sigma_S(z)= \frac{\|K\|_2^2\hat \sigma_S^2(z)}{\hat f(z)}   ~~ \mbox{and}~~ \hat\Sigma^*_S(z)= \frac{\|K\|_2^2 \hat \sigma_S^{*2}(z)}{\hat f(z)},
 \end{align}
 where $\hat \sigma_S(z)$ and $\hat \sigma_S^{*2}$ are estimators for $ \sigma_S(z)$ and $\sigma_S^{*2}$ with
\begin{align*}
 & \hat \sigma_S^2=\frac{1}{nh^l}\sum_{i=1}^n \frac{  (\psi(\hat q, \mathcal{W}_i)-\hat \tau_S(z))^2K_i} {\hat f(z)}~\mbox{and}\\
 &  \hat \sigma_S^{*2}=\frac{1}{nh^l}\sum_{i=1}^n \frac{  (\psi^*(\hat q, \mathcal{W}_i)-\hat \tau_S(z))^2K_i} {\hat f(z)},
 \end{align*}
 $\hat f(z)= {\sum_{i=1}^nK_h(Z_i-z)}/n$ is a  kernel-based estimator of $f(z)$,
 \begin{align*}
 &\psi(\hat q,\mathcal{W}_i) = \frac{D_iY_i}{\hat q(V^{\top}X_i)}
    - \frac{(1-D_i)Y_i}{ 1-\hat q(V^{\top}X_i) }   ~~ \mbox{and}~~
 \psi^*(\hat q,\mathcal{W}_i) =  \frac{[D_i\{Y_i-\hat m_1(V^{\top}
    X_i)\} ]}{ \hat q(V^{\top} X_i)} \\
    &- \frac{[(1-D_i)\{Y_i-\hat m_0(V^{\top} X_i)\} ]}{1- \hat q(V^{\top} X_i)}
    +\hat m_1(V^{\top} X_i)-\hat m_0(V^{\top} X_i)
 \end{align*}
with  $\hat m_1(\hat V^{\top}X)= {\sum_{\{t: D_t=1\}}^n\mathcal{H}_{h_2}(\hat V^{\top}X_t-\hat V^{\top}X)Y_t}\Big/{\sum_{\{t: D_t=1\}}^n\mathcal{H}_{h_2}(\hat V^{\top}X_t-\hat V^{\top}X)}$,  $\hat m_0(\hat V^{\top}X)= {\sum_{\{t: D_t=0\}}^n\mathcal{H}_{h_2}(\hat V^{\top}X_t-\hat V^{\top}X)Y_t}\Big/{\sum_{\{t: D_t=0\}}^n\mathcal{H}_{h_2}(\hat V^{\top}X_t-\hat V^{\top}X)}$ being the estimators of $m_1(V^{\top}X)$ and $m_0(V^{\top}X)$.

Further we can state the consistency of the proposed  estimators in the following theorem.
\begin{theorem}
Suppose all the conditions in Appendix are satisfied, we have that
\begin{align*}
 \hat \Sigma_S(z)= \Sigma_S(z)+o_p(1),   ~~ \mbox{and}~~ \hat\Sigma^*_S(z)=\Sigma_S^*(z) +o_p(1).
\end{align*}
\end{theorem}

By Theorem~2, we can obtain the pointwise consistent estimator for standard error of $\sqrt{nh^l}(\hat \tau_S(z)-\tau(z))$, so that we are able to construct a $(1-\alpha)100\%$ pointwise confidence interval for $\tau(z)$, i.e.
\begin{align}\label{cib_1}
\hat \tau_S(z)\pm (nh)^{-1/2}c_{\alpha/2}\left(\hat \Sigma_S(z)\right)^{1/2}
\end{align}
or
\begin{align}\label{cib_2}
\hat \tau_S(z)\pm (nh)^{-1/2}c_{\alpha/2}\left(\hat \Sigma_S^*(z)\right)^{1/2},
 \end{align}
 with $c_{\alpha/2}$ being the $(1-\alpha/2)$ quantile of the standard normal distribution. Note that the specification formula of confidence interval depends on whether the condition $|Z\cap V^{\top}X|<l$ or $|Z\cap V^{\top}X|=l$.
 One possible way to make choice between \eqref{cib_1} and \eqref{cib_2} is based on the value of $|Z\cap \hat V^{\top}X|$.

 To be specified, taking MAVE\citep{xia2002} as an exmaple dimension reduction method, we proposed a estimation and inference
procedure of $\tau(z)$ based on $IPW$-$S$ by carrying out the following steps.

\begin{description}
\item {Step 1:} Obtain the estimator of $V$  by solving the minimizing problem
\begin{eqnarray*}
\min\limits_{V,a,b}\sum_{i,j=1}^n\{D_i-a_j-b_j^{\top}V^{\top}(X_i-X_j)\}^2\omega_{ij}.
\end{eqnarray*}
Here $\omega_{ij}=\mathcal{H}_{h_2}\{V^{\top}(X_i-X_j)\}\big /\sum_{l=1}^n\mathcal{H}_{h_2}\{V^{\top}(X_l-X_j)\}, a=(a_1,\ldots,a_n)$, $b=(b_1,\ldots,b_n)$.
Denote the resulting estimator by $\hat V$

\item {Step 2:} Given $\hat V$, estimate the propensity score $E(D\mid \hat V^{\top}X)$ via \eqref{semi_q}.

\item {Step 3:} Obtain the semiparametric CATE estimator, $\hat \tau_S(z)$, via \eqref{est_scate}.

\item {Step 4:} Given $Z=z$, a $(1-\alpha)$ pointwise confidence interval for the true CATE, $\tau(z)$, can be constructed as follows.
  if $|Z\cap \hat V^{\top}X|\approx l$, the confidence interval of $\tau(z)$ can be constructed in the form of \eqref{cib_2}, i.e. $
\hat \tau_S(z)\pm (nh)^{-1/2}c_{\alpha/2}\left(\hat \Sigma_S^*(z)\right)^{1/2},$
   Otherwise, the pointwise confidence interval of $\tau(z)$ would be constructed in the form of \eqref{cib_1}, that is, $
\hat \tau_S(z)\pm (nh)^{-1/2}c_{\alpha/2}\left(\hat \Sigma_S(z)\right)^{1/2}.$
\end{description}
Note that the first step can be implemented using the R package MAVE. Based on this estimation and inference procedure, the empirical analysis in section~5 can be implemented.
 }

\subsection{ Extension of existing IPW-N}

 Recall $IPW$-$N$ proposed by \cite{Abrevaya2015} is
\begin{eqnarray}\label{est_ncate}
\hat \tau_N(z)=\left( { \sum_{i=1}^n
\left[\frac{D_iY_i}{\hat p( X_i)} -\frac{(1-D_i)Y_i}{1-\hat
 p (X_i)}\right]K_{h}(Z_i-z)}\right)\Big /
{ \sum_{i=1}^nK_{h}(Z_i-z)},
\end{eqnarray}
 with
 $\hat p( X_i)= { \sum_{j\neq i}^nD_j
\mathcal{L}_{h_1}({X_j-  X_i})}  \Big /
\sum_{j\neq i}^n\mathcal{L}_{h_1}({X_j-X_i}).$ {Here $\mathcal{L}(\cdot)$ is also a multivariate kernel function with $\mathcal{L}_{h_1}(\cdot)=h_1^{-k}\mathcal{L}(\cdot/h_1)$, and $h_1$ is the corresponding bandwidth.}

Note that the asymptotic properties of $IPW$-$S$ is influenced by the affiliation of $Z$, we  in this paper try to analyse the asymptotic properties of $IPW$-$N$ in different  scenarios similarly as  the ones in Theorem \ref{sem_them}. Suppose
  \begin{eqnarray}\label{non_assum}
  D\perp X\mid \tilde X, ~\tilde X
\subseteq X,~\tilde k=dim(\tilde X)\leq k.
\end{eqnarray}   To extend the  asymptotic results of $IPW$-$N$ in \cite{Abrevaya2015}, we derive the following theorem that also confirms the influence of the affiliation of $Z$ to $\tilde X$ in the asymptotic properties of $IPW$-$N$.  \cite{Abrevaya2015} only considered a special situation in the following Theorem~\ref{nonp_them}: $|Z \cap \tilde X|=l$ and $\tilde X=X$.

Before stating the result as theorem, let us define some important quantities:
\begin{align*}
&\psi(p(\tilde
X),\mathcal{W})= \frac{DY}{p(\tilde X)}
    -\frac{(1-D)Y}{1-p(\tilde X)} ,\\
&\psi^*(p(\tilde
X),\mathcal{W})=\frac{[D\{Y-m_1(\tilde
    X)\}]}{ p(\tilde X)}-\frac{[(1-D)\{Y-m_0(\tilde X)\}]}{1- p(\tilde X)}
    +m_1(\tilde X)-m_0(\tilde X).
\end{align*}

\begin{theorem}{\label{nonp_them}}
Suppose all the conditions in Appendix are satisfied, the
following statements  hold for each point $z$ in the support of $Z$:
\begin{itemize}
 \item[(1)] { When $|Z \cap \tilde X|=t<l$ with
    $s_1[2-l/(l-t)]+l>0$, the asymptotically linear
    representation is
 \begin{eqnarray*}
\sqrt{nh^l}\left(\hat \tau_N(z)-\tau(z)\right)= \frac{1}{\sqrt{nh^l}f(z)} \sum_{i=1}^{n}[\psi(p(\tilde X_i),\mathcal{W}_i)-\tau(z)]K_i+o_p(1);
\end{eqnarray*}
 the asymptotic distribution of $\hat \tau_N(z)$  is
\begin{eqnarray*}
\sqrt{nh^l}\left(\hat \tau_N(z)-\tau(z)\right)
\stackrel{\mathcal{D}}{\longrightarrow}N\left(0,\Sigma_N(z)
 \right).
\end{eqnarray*}}
 \item[(2)]When $|Z \cap \tilde X|=l$,  the asymptotically linear
    representation is
 \begin{eqnarray*}
\sqrt{nh^l} \left(\hat \tau_N(z)-\tau(z)\right)= \frac{1}{\sqrt{nh^l}f(z)} \sum_{i=1}^{n}[\psi^*(p(\tilde X_i),\mathcal{W}_i)-\tau(z)]K_i+o_p(1);
\end{eqnarray*}
 the asymptotic distribution of $\hat \tau_N(z)$ is
\begin{eqnarray*}
\sqrt{nh^l} \left(\hat \tau_N(z)-\tau(z)\right)
\stackrel{\mathcal{D}}{\longrightarrow}N\left(0, \Sigma_N^*(z)
\right).
\end{eqnarray*}
\end{itemize}
Here $s_1$ is the order of $\mathcal{L}(\cdot)$,  $\Sigma_N(z)= {\parallel
K \parallel_2^2\sigma_{N}^2(z)}\Big /{f(z)}$, $\Sigma^*_N(z)={\parallel
K \parallel_2^2\sigma_{N}^{*2}(z)}\Big/{f(z)}$,
 $\sigma_{N}^{2}(z)
    =E[\{\psi(p(\tilde X_i),\mathcal{W}_i)-\tau(z)\}^2\mid
Z=z]$, $\sigma_{N}^{*2}(z)=E[\{\psi^*(p(\tilde
X_i),\mathcal{W}_i)-\tau(z)\}^2\mid Z=z]$.

\end{theorem}

{
Similarly as $IPW$-$S$, we also proposed the estimators for $\Sigma_N(z)$ and $\Sigma^*_N(z)$  under $|Z\cap \tilde X|<l$ and $|Z\cap \tilde X|=l$ respectively as
 \begin{align}
 \hat \Sigma_N(z)= \frac{\|K\|_2^2\hat \sigma_N^2(z)}{\hat f(z)},   ~~ \mbox{and}~~ \hat\Sigma^*_N(z)= \frac{\|K\|_2^2 \hat \sigma_N^{*2}(z)}{\hat f(z)},
 \end{align}
 where $\hat \sigma^2_N$ and $\hat \sigma_N^{*2}$ are estimators for $ \sigma^2_N(z)$ and $\sigma_N^{*2}(z)$ with
\begin{align*}
 & \sigma_N^2=\frac{1}{nh^l}\sum_{i=1}^n \frac{ ( \psi(\hat p, \mathcal{W}_i)-\hat \tau_N(z))^2K_i} {\hat f(z)}~,~ \hat \sigma_N^{*2}=\frac{1}{nh^l}\sum_{i=1}^n \frac{ ( \psi^*(\hat p, \mathcal{W}_i)-\hat \tau_N(z))^2K_i} {\hat f(z)},\\
 &\psi(\hat p,\mathcal{W}_i) = \frac{D_iY_i}{\hat p(\tilde X_i)}
    -  \frac{(1-D_i)Y_i}{ 1-\hat p(\tilde X_i) } ~\mbox{and}~~
\psi^*(\hat p,\mathcal{W}_i)
 = \frac{[D_i\{Y_i-\hat m_1(\tilde
    X_i)\} ]}{ \hat p(\tilde X_i)}
    \\&  - \frac{[(1-D_i)\{Y_i-\hat m_0(\tilde X_i)\} ]}{1- \hat p(\tilde  X_i)}
    +\hat m_1(\tilde X_i)-\hat m_0(\tilde X_i).
 \end{align*}
And $ \hat m_1(\tilde X)={\sum_{\{t: D_t=1\}}^n\mathcal{L}_{h_1}(\tilde X_t-\tilde X)Y_t}/ {\sum_{\{t: D_t=1\}}^n\mathcal{L}_{h_1}(\tilde X_t-\tilde X)}$,
  $\hat m_0(\tilde X)=\\{\sum_{\{t: D_t=0\}}^n\mathcal{L}_{h_1}(\tilde X_t-\tilde X)Y_t} /{\sum_{\{t: D_t=0\}}^n\mathcal{L}_{h_1}(\tilde X_t-\tilde X)}$ {being the estimators of}~ $m_1(\tilde X)$ ~\mbox{and}~$m_0(\tilde X)$.
Further we can show the consistency of proposed asymptotic variance function estimators via the following theorem.
\begin{theorem}
Suppose all the conditions in Appendix are satisfied, we have that
\begin{align*}
 \hat \Sigma_N(z)= \Sigma_N(z)+o_p(1),   ~~ \mbox{and}~~ \hat\Sigma^*_N(z)=\Sigma_N^*(z) +o_p(1).
\end{align*}
\end{theorem}
\begin{remark}
Based on Theorem~4, we can also get the consistent estimator for standard error of $\sqrt{nh^l}(\hat \tau_N(z)-\tau(z))$ and  construct a
pointwise confidence interval of $\tau(z)$ based on $\hat \tau_N(z)$. However, we first need to estimate the true active arguments of propensity score $\tilde X$, denoting the corresponding estimator as ${\hat X}$, which can be done by variable selection method, to decide the proper form of confidence interval. To be specified, if $|Z\cap \hat X|\approx l$, the pointwise confidence interval can be constructed as
\begin{align*}
\hat \tau_N(z)\pm (nh)^{-1/2}c_{\alpha/2}\left(\hat \Sigma_N^*(z)\right)^{1/2}.
\end{align*}
Otherwise, we would construct the pointwise confidence interval as
\begin{align*}
\hat \tau_N(z)\pm (nh)^{-1/2}c_{\alpha/2}\left(\hat \Sigma_N(z)\right)^{1/2}.
 \end{align*}

\end{remark}
}

\subsection { Some further studies on estimation efficiency }
When $\tilde X=X$, as proved by \cite{Abrevaya2015}, $IPW$-$N$ can be asymptotically more efficient than $IPW$-$P$:
\begin{eqnarray*}
&&\sigma_{P}^2(z)=\sigma_{N}^{2*}(z)+
E\left[p(X)\{1-p(X)\}\left\{\frac{m_1(X)}{p(X)}
+\frac{m_0(X)}{1-p(X)}\right\}^2 \mid Z=z\right],
\end{eqnarray*}
and $IPW\mbox{-}P\cong IPW\mbox{-}O.$
Here $m_j(X)=E\{Y(j)\mid X\}.$
Thus, with $p(X)=p(\tilde X)=q(V^{\top}X),$ we can give the ranking for the estimation
  efficiency of the four estimators in the following corollary.
\begin{corollary}\label{cor_sem}Suppose all the assumptions and conditions in Appendix are satisfied and $p(X)=p(\tilde X)=q(V^{\top}X)$, the
following statements  hold for each point $z$ in the support of $Z$:

 {\bf{Case~1:}}  When $|Z \cap \tilde X|=l$ with $\tilde X=X$ and $|Z \cap V^{\top}X|=l$, $$IPW\mbox{-}N \preceq IPW\mbox{-}S \preceq IPW\mbox{-}P \cong IPW\mbox{-}O,~\mbox{with}$$
\begin{eqnarray*}
\sigma_{P}^2(z)=\sigma_{S}^{*2}(z)+
E\bigg[q(V^{\top}X)(1-q(V^{\top}X))
\bigg\{\frac{m_1(V^{\top}X)}{q(V^{\top}X)}
 +\frac{m_0(V^{\top}X)}{1-q(V^{\top}X)}\bigg\}^2
\mid Z=z\bigg],\\
\sigma_{S}^{*2}(z)=\sigma_{N}^{*2}(z)+
E\bigg[q(V^{\top}X)(1-q(V^{\top}X))\bigg\{\frac{\Delta m_1}{q(V^{\top}X)}
 +\frac{\Delta m_0}{1-q(V^{\top}X)}\bigg\}^2\mid Z=z\bigg],\quad
\end{eqnarray*}
where $\Delta m_j=m_j(X)-m_j(V^{\top}X).$

 {\bf{Case~2:}} When $|Z \cap \tilde X|=l$ with $\tilde X=X$ but $|Z \cap
    V^{\top}X|=t$ with $0\leq t <l$,  $$IPW\mbox{-}N \preceq IPW\mbox{-}S \cong IPW\mbox{-}P \cong IPW\mbox{-}O$$ with
    $\sigma_{S}^2(z)=\sigma_{P}^2(z)=\sigma_{O}^2(z).$

 {\bf{Case~3:}} When $|Z \cap \tilde X|=t$ with $\tilde X \subsetneq X$ and $|Z
    \cap V^{\top}X|=t$ with $0\leq t < l$,   $$IPW\mbox{-}N \cong IPW\mbox{-}S \cong IPW\mbox{-}P \cong
    IPW\mbox{-}O$$ with
    $\sigma_{N}^2(z)=\sigma_{S }^2(z)=\sigma_{P}^2(z)=\sigma_{O }^2(z).$

\end{corollary}
\begin{remark}
 In {\bf Case 1}, the
equality in the first inequality holds when both $m_1(V^{\top}X)$ and
$m_0(V^{\top}X)$ equal to zero, and the equality in the second
 inequality holds when $m_j(X)=m_j(V^{\top}X)$ for all $j=0,1$. A sufficient condition to make $m_j(X)=m_j(V^{\top}X)$ hold is $E(Y_j\mid
X)\perp X\mid V^{\top}X$ meaning that $Y(1)$ and $Y(0)$ share the same
central mean subspace.
\end{remark}

\begin{remark}
Here, we  discuss another special case: $V^{\top}X=Z$  in Corollary~\ref{cor_sem} such that $q(V^{\top}X)=p(Z)$. It follows
that $IPW\mbox{-}S \preceq IPW\mbox{-}P$ with $
 \sigma_{P}^2(z)=\sigma_{S}^{*2}(z)+
p(z)(1-p(z))
\left[m_1(z)/\{p(z)\}
+m_0(z)/\{1-p(z)\} \right]^2.$
 Similarly, $IPW\mbox{-}N \preceq IPW\mbox{-}P$ if
$ \tilde X=Z$: $
 \sigma_{P}^2(z)=\sigma_{N}^{*2}(z)+
p(z)(1-p(z))\left[m_1(z)/\{p(z)\}
+m_0(z)/\{1-p(z)\} \right]^2.$
Thus, if $Z=V^{\top}X=\tilde X$, we have $
\sigma_{S}^{*2}(z)=\sigma_{N}^{*2}(z)\leq \sigma_{P}^2(z).$
\end{remark}

\begin{remark}
  Although $IPW$-$S$ cannot be more efficient than $IPW$-$N$ in theory, it
has an obvious advantage due to its dimension reduction structure. This can be very useful
in practice as when $X$ is of high dimension, $IPW$-$N$ is hard to use as it
has to adopt very high order kernel function and delicately chosen
bandwidths. {The numerical studies in the next section show that when the
dimension of $X$ is only $4$, $IPW$-$S$ can performs better than $IPW$-$N$ in some cases. Thus, in the numerical studies, when the in high
dimension $k=20$, we do not consider $IPW$-$N$.}

{Another issue is also relevant. Generally speaking, combining the results in Subsections~3.1 and 3.3, when the dimension reduced vector $V^{\top}X$ cannot fully cover the given covariates $Z$, the $IPW$-$S$ is less efficient. It seems that we can add $Z$ into the covariates to be $(Z, V^{\top}X)$ to enhance the estimation efficiency in theory. However, this causes the estimation procedure much more complicated (with higher order kernel and more delicately selected bandwidths) and less accurate due to the dimension increasing as described above. Thus, balancing the theoretical merit and practical usefulness, we still prefer using $IPW$-$S$ without adding more covariates.}
\end{remark}

{
From the above discussion, we can find that the asymptotic efficiency comparison result of IPW-type CATE estimators is different from that of IPW-type ATE estimators, because
$\mbox{ATE estimator using nonparametric esitmated}~p(X) \preceq \mbox{that using parametric} \\\mbox{estimated}~p(X) \preceq \mbox{that using true}~p(X)$.}
Thus it is worthwhile to give a further exploration on the reasons. {From our study, we find that it is mainly because of the different convergence rates of the estimated propensity scores under different
scenarios.}
In the following corollary, we show that when the convergence rate of the nonparametically estimated
propensity score can be fast enough, $IPW$-$N$ and $IPW$-$S$ can also be
asymptotically equivalent to $IPW$-$P$, so is $IPW$-$O$. This is the case when the
propensity score function is smooth sufficiently and the
 kernel and
bandwidths are chosen delicately to meet the mentioned condition in Corollary~\ref{cor2}.

\begin{corollary}\label{cor2}
Suppose all the conditions in Appendix are satisfied.
\begin{itemize}
  \item[(1)] When
$\sqrt{nh^l} \left(h_2^{s_2}+\sqrt{\log(n)/nh_2^r}\right)=o(1)$, it follows that $$ IPW\mbox{-}S \cong IPW\mbox{-}P \cong IPW\mbox{-}O. $$
  \item[(2)] When
$\sqrt{nh^l} \left(h_1^{s_1}+\sqrt{\log(n)/nh_1^{\tilde k}}\right)=o(1)$, it follows that $$IPW\mbox{-}N  \cong IPW\mbox{-}P \cong IPW\mbox{-}O. $$
  \item[(3)] When
 $\sqrt{nh^l} \left(h_2^{s_2}+\sqrt{\log(n)/nh_2^r}\right)=o(1)$ and $\sqrt{nh^l} \left(h_1^{s_1}+\sqrt{\log(n)/nh_1^{\tilde k}}\right)=o(1)$, it follows that $$IPW\mbox{-}N \cong IPW\mbox{-}S \cong IPW\mbox{-}P \cong IPW\mbox{-}O. $$
\end{itemize}
\end{corollary}

\begin{remark}
{Corollary~\ref{cor2} implies that when the convergence rate of estimated propensity score is fast enough, the corresponding CATE estimator would be asymptotically equivalent to $IPW$-$O$, which is based on true propensity score, even though the condition $|Z \cap V^{\top}X|$ in Theorem~\ref{sem_them} or $|Z \cap \tilde X|$ in Theorem~\ref{nonp_them} is satisfied. In this sense, we can say that the convergence rate of estimated propensity score is dominant the role of affiliation of $Z$ in the set of arguments of propensity score in comparing the asymptotic efficiencies among the CATE estimators.}
It is well known that the convergence rate of nonparametric estimator is
possibly close to ${n}^{-1/2}$  if the estimated function is very smooth
and the higher kernel function is utilized, see \cite{li2007}. Thus the
conditions $\sqrt{nh^l} \left(h_2^{s_2}+\sqrt{\log(n)/nh_2^r}\right)=o(1)$ and $\sqrt{nh^l} \left(h_1^{s_1}+\sqrt{\log(n)/nh_1^{\tilde k}}\right)=o(1)$ could hold.
As the choices for such kernel and bandwidths often make no sense
for practical use, this investigation only serves as a theoretical
exploration with a remind that a ``good" estimator for the propensity score
may not be helpful for constructing a ``good" CATE.

\end{remark}

\section{ Simulation study}

 \subsection { Preliminary of the simulation}
To evaluate the finite sample performances of
$IPW$-$S$, we consider the comparisons with
 $IPW$-$P$, $IPW$-$N$ and $IPW\mbox{-}O.$  {
To save space, we only present the simulations in the case $Z\in R$. To make the comparisons more convincing, we consider
two scenarios with two low dimensions of $X=(Z, U_1,\cdots,U_{k-1})$ equal to $k=2$ and $4$, and higher dimensions $k=20$. In the latter,
   $IPW$-$N$ is not included as very high order kernel and very delicately
    selected bandwidths are required and then it is very difficult to implement.
Several criteria are used to evaluate the estimation
  efficiency: $Bias$; estimated standard deviation $Est\_SD$; mean
  square error (MSE). As the asymptotic distributions are
   standard normal, we then also report the proportions
   outside
 the critical values $\pm 1.645$: $P_{\pm 1.645}$.
Further, to make the efficiency ranking in finite sample setting more
visible, we   report, as relative efficiency, the $Est\_SD$ results via dividing each $Est\_SD$
by $Est\_SD$ of $IPW$-$O$ that is used as the benchmark.
Thus, when the ratio is smaller than $1$, the
corresponding estimator is more efficient than $IPW$-$O$.

\subsection { Experiment 1(low-dimensional setting)}

In the low-dimensional setting, we consider the covariates $X=(Z, U_1,\cdots,U_{k-1})$ are given by the following procedure.
When $k=2$, $X=(Z,U_1)$ with $
Z=\epsilon_1$ and $U_1=(1+2Z)^2(-1+Z)^2+\epsilon_2.
$
When $k=4$, $X=(Z,U_1,U_2,U_3)$ are given by $
  Z=\epsilon_1,~ U_1=(1+2Z)+\epsilon_1,~
  U_2=(1+2Z)+\epsilon_2,~ U_3=(-1+Z)^2+\epsilon_3.$ $\epsilon_i\sim  unif[-0.5,0.5]$ for $i=1,2,3,$ and they are mutually independent.
 To easily compare the theoretical results under parametric, nonparametric and semiparametric structure,
 we consider four models:
 {\begin{itemize}
 \item {{Model 1 (k=2, r=1 with $|Z\cap  X|=1$ but $|Z\cap V^{\top}X|=0$):} } \begin{eqnarray*}
  Y(1) =\beta_1^{\top}X+\gamma_1ZU_1+\nu,~Y(0)=0 ~ \mbox{and}~ p_1(X)=\Lambda(1/\sqrt{2}(Z+U_1)).\end{eqnarray*}

 \item{{Model 2 (k=2, r=1 with $|Z\cap \tilde X|=|Z\cap U_1|=0$ and $|Z\cap V^{\top}X|=0$):}} \begin{eqnarray*}
  Y(1)=\beta_1^{\top}X+\gamma_1ZU_1+\nu,~Y(0)=0~\mbox{and}~ p_2(X)=\Lambda(U_1).
  \end{eqnarray*}

 \item {{Model 3 (k=4, r=1 with $|Z\cap X|=1$ but $|Z\cap V^{\top}X|=0$):}}
    \begin{eqnarray*}
 Y(1)=\beta_2^{\top}X+\gamma_2ZU_1U_2U_3+\nu, ~Y(0)=0, ~\mbox{and} ~p_3(X)=\Lambda(0.5(Z+U_1+U_2+U_3)).
\end{eqnarray*}

 \item {{Model 4 (k=4, r=2 with $|Z\cap X|=1$ and $|Z\cap V^{\top}X|=1.$):}}\begin{eqnarray*}
 Y(1)=\beta_2^{\top}X+\gamma_2 ZU_1U_2U_3+\nu,~Y(0)=0~\mbox{and}~
p_4(X)=\Lambda\left(\frac{\sqrt{3}(1+Z)}{\sqrt{3}+U_1+U_2+U_3}\right).
\end{eqnarray*}
\end{itemize}
}

Here
 $\nu \sim N(0,0.25^2)$,
$\Lambda(\cdot)$ is the c.d.f. of the logistic distribution. Given that the matrix $V$ satisfies $E(D|X)\perp X|V^{\top}X$, we consider these four types of propensity score models to satisfy the
 conditions in different scenarios. Under Model~1 and Model~3,  the
 dimension of $V^{\top}=(1,\cdots,1)_{1\times k}$ is $dim(V)=r=1$, 
$|Z\cap  X|=1$ but $|Z\cap V^{\top}X|=0$.
  Thus, we aim to examine
 whether $IPW\mbox{-}N \preceq IPW\mbox{-}S \cong IPW\mbox{-}P \cong IPW\mbox{-}O$. {In order to examine the theoretical
 results
in {\bf Case 3} of Corollary~\ref{cor_sem}  in the finite sample
 scenario, we  consider $ p_2(X)$ in Model 2. In this setting,
  $D\perp X |\tilde X=U_1$, 
and for $V^{\top}=(0,1)$,  $ p_2( X)\perp X
|V^{\top}X$. Obviously, $|Z\cap \tilde X|=|Z\cap U_1|=0$ and $|Z\cap V^{\top}X|=0$ in Model~2, thus it can be used to examine whether $IPW\mbox{-}N  \cong IPW\mbox{-}S\cong IPW\mbox{-}P  \cong IPW\mbox{-}O $.
$ p_4(X)$ in Model~4 is also set to verify the results
in Corollary \ref{cor_sem}. This propensity score function has
   $Z$ itself as an individual argument. Namely, $p_4(X)\perp X\mid V^{\top}X$
    with  $V^{\top}=\left(
                 \begin{array}{cccc}
                   1 & 0 & 0 & 0\\
                   0 & 1/\sqrt{3}& 1/\sqrt{3}& 1/\sqrt{3} \\
                 \end{array}
               \right)$ and $dim(V)=r=2$, and $|Z\cap X|=1$ and $|Z\cap V^{\top}X|=1.$ We will examine whether
  $IPW$-$S$ can be  more
   efficient than $IPW$-$P$ and $IPW$-$O$.
   {
 As for the parameters $\beta_i$, $\gamma_i$, $i=1,2$, we consider following scenarios:
 \begin{itemize}
 \item Scenario I:
 $\beta_1^{\top}=(0,0)$, $\gamma_1=1$, $\beta_2^{\top}=(1/10,1/\sqrt{2},-1/\sqrt{2},-1/10)$, $\gamma_2=0$;
 \item Scenario II:
 $\beta_1^{\top}=(1/2,-1/5)$, $\gamma_1=0$, $\beta_2^{\top}=(0,0,0,0)$, $\gamma_2=1$.
 \end{itemize}
 Obviously, when $r_i=0$, the linear model is being considered, while $r_i=1$, the nonlinear model is taken into account, $i=1,2$.
   }

   Next,  we  determine the order of kernels $\mathcal{L}(\cdot)$, $\mathcal{H}(\cdot)$
   and $K(\cdot)$ to guarantee the regularity
   condition~5 in Appendix.
  As
  there is no data-driven or optimal selection method available for $IPW$-$N$ and $IPW$-$S$, we use the  rule of
  thumb to select them as  suggested by \cite{Abrevaya2015} for fair comparisons.
   The principle of selection is based on  proper
    rates of convergence  in the form of $h=a\cdot n^{-\eta}$ for $a>0$
    and $\eta >0$.
    Since $IPW$-$S$ can be regarded as a low-dimensional
   type $IPW$-$N$, the bandwidths can be chosen via replacing $\tilde k$ by $r$ as follows:
   \begin{eqnarray}\label{sem_band_sele}
 h_2=a_2 n^{ {-1}/(2r+\delta_r+\delta_2)},  \
  h=a n^{ {-1}/(l+4+2r+2\delta_r-\delta_1)},
   \end{eqnarray}

where $a,a_2,\delta_2,\delta_1$ are positive. Note that $\delta_2$ and $\delta_1$ can be as small as desired,  thus we let them
to be zero in the simulations for simplicity. Further, the order of
$\mathcal{H}$ is $s_2=r+\delta_r$ with $\delta_r=0$ for even $r$, $\delta_r=1$
for odd $r$, and $s=s_2+2=r+\delta_r+2$. Due to the semiparametric nature, we
construct two consistent estimators  $\hat V$ via  MAVE proposed by \citep{xia2002, xia2007}, which can be implemented in the R package MAVE.

 Further, to
fairly  examine the performances, the parameters $s,h$ for $K(u)$
are the same for all these four $CATE$ estimators. Since $r \leq \tilde k$, the
choices of $s$ and $h$ for $IPW$-$N$ can be used for  all the $CATE$
estimators.
Taking all into account, the corresponding bandwidths
are
summarized in Table~\ref{turnpara}.
\begin{table}[H]
\caption{The order of bandwidths in the simulations.}
 \label{turnpara}\par
   \renewcommand\arraystretch{0.8}
\vskip .2cm
\centerline{\tabcolsep=10truept\begin{tabular}{cccc} \hline
  $\tilde k$ &  $h$ & $h_1$ & $h_2$\\
  \hline
   $\tilde k$=1 & $a\cdot n^{-1/9}$ & $a_1\cdot n^{-1/3}$ & $a_2\cdot n^{-1/(2 r+\delta_r)}$  \\
  $\tilde k$=2& $a\cdot n^{-1/9}$ & $a_1\cdot n^{-1/4}$ & $a_2\cdot n^{-1/(2  r+\delta_r)}$  \\
 $\tilde k$=4& $a\cdot n^{-1/13}$ & $a_1\cdot n^{-1/8}$ & $a_2\cdot n^{-1/(2  r+\delta_r)}$  \\
   \hline
\end{tabular}}
\end{table}
 As for the tuning parameters ${a,a_1,a_2}$, we consider the following two groups of  values: $Group\ 1: \{a=0.55,a_1=1.05,a_2=0.75\}$, $Group\ 2:~
\{a=0.55,a_1=1.05,a_2=0.69\}.$

 Specifically,  we estimate $\tau(Z) $ at  $Z\in\{-0.4,-0.2,0,0.2,0.4\}$. The
sample sizes are $n=500$ and $1000.$ The replication time is $500$. We
choose the Gaussian kernel  and higher order kernels derived from it
throughout this section. Further, we should point out that the estimated
propensity score is trimmed to lie in the interval $[0.005,0.995]$ as \cite{Abrevaya2015} did.
{We give the observations from the simulation results reported in
Tables~\ref{k=2_r=1_G1_s1}-\ref{k=4_r=2_s2}.
To save space, we only report the results about the the relative efficiency with $\{a=0.55,a_1=1.05,a_2=0.75\}$ under Scenario I. See Figure~$\ref{k=low}$. We have the following observations.}

{\it Observation 1.}  As  expected, larger sample size leads to smaller bias and standard
deviation in most cases. When $k=4$ and the sample size goes
from $n=500$ up to $n=1000$, the
bias and variance reduction are more significant and the empirical values of $P_{1.645}$
and $P_{-1.645}$ are closer to the nominal level $0.05$. That
implies that the normal approximation works well.

{\it Observation 2.}  When the dimension $k$ increases  to $4$,
$Bias $, $Est\_SD$ and $MSE$ also increase. Further, the  dimension
does have impact on the performance of $IPW$-$N$. {In Tables~$S.1$} under
Model~1 with $k=2$, $IPW$-$N$ is uniformly more efficient than all the others. But under the models with $k=4$, especially  with $k=4$ and $r=2$,
the superiority of $IPW$-$N$ becomes less significant. We can
see that  $IPW$-$S$ can even be more efficient than $IPW$-$N$ sometimes, mainly
due to its dimension reduction structure.

{

{\it Observation 3.} Taking into account all the simulation results in Tables \ref{k=2_r=1_G1_s1}-\ref{k=4_r=2_s2}., the estimated standard deviation ($Est\_SD$) of all CATE estimations increase as $z$ is close to the boundary of the support of $Z$.  This phenomenon should be  mainly because of  the nonparametric estimation of CATE(z) function with respect to $Z$. Note that $IPW$-$O$ also involves nonparametric estimation for the conditional expectation over $Z$, thus the boundary effect also takes place for it.
 Further, empirically,  the $Est\_SD$ of $IPW$-$O$ often increases, in the numerical studies we conduct, relatively more quickly than $IPW$-$P$  or $IPW$-$S$ in the cases we will discuss in Observation 4 below when $z$ is close to the boundary. Figure~\ref{k=low} about the relative efficiency compared with $IPW$-$O$ shows this though for different models, at the boundary,  $IPW$-$O$ has different relative efficiency with $IPW$-$P$. Combining the  information that ATE with estimated unknowns, the final estimators could be more efficient in general, less relative efficiency of $IPW$-$O$ in finite sample scenarios may be understandable although in the case that the asymptotic efficiency should be equivalent. On the other hand, when we look at the original values of  $IPW$-$O$, the differences with $IPW$-$P$ is not significant.  

{\it Observation 4}. We also check the effect caused by the inclusiveness of the given
covariates $Z$ in the set of the arguments of the propensity score for
$IPW$-$S$ and $IPW$-$N$.    Under Model~1 or Model~3 with  $|Z\cap  \tilde X|=1$ but $|Z\cap V^{\top} X|=0$, Figure $\ref{k=low}$
 shows that  the $Est\_SD$
 of $IPW$-$N$ is uniformly smaller than those of the other $CATE$ estimators.
 While, $IPW$-$S$ and $IPW$-$P$ have similar performance. This coincides with the theory.
In contrast, under Model~2 where $|Z\cap \tilde X|=|Z\cap U_1|=0$, $IPW$-$N$ losses its superiority of efficiency to share similar performance to $IPW$-$S$ and $IPW$-$P$. 
Under Model~4 where $|Z\cap X|=1$ and $|Z\cap V^{\top}X|=1$,  $IPW$-$S$  outperforms $IPW$-$O$ and $IPW$-$P$, and can even be comparable with $IPW$-$N$ sometimes. These
results also coincide with the theory in Corollary \ref{cor_sem}.

}


%
%

\begin{figure}[H]
\flushleft
  \includegraphics[scale=0.43]{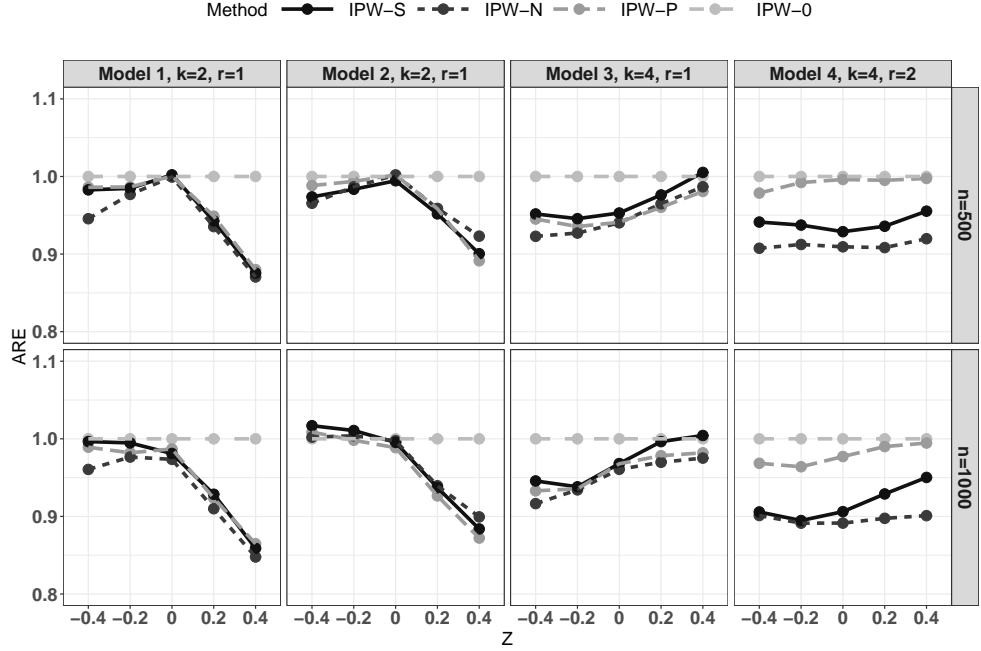}\\
\caption{The asymptotic relative efficiency(ARE) about $Est\_SD$ against that of $IPW$-$O$ under $Group~1: \{a=0.55,a_1=1.05,a_2=0.75\}$ and Scenario I:
 $\beta_1^{\top}=(0,0)$, $\gamma_1=1$, $\beta_2^{\top}=(1/10,1/\sqrt{2},-1/\sqrt{2},-1/10)$, $\gamma_2=0$.
}\label{k=low}
\end{figure}

\subsection { Experiment 2 (high dimensional  setting).}

Consider models with  much higher dimensional $X$:
$k=dim(X)=20$. As $IPW$-$N$ obviously suffers  from the curse of
dimensionality and thus does not work at all, we then only focus on $IPW$-$O$,
$IPW$-$P$ and $IPW$-$S$. 
To better examine the corresponding finite sample performances, we consider
the model settings which are similar to  Models~3 and~4  with uniformed $Z$, but  with more zero
coefficients for ease of comparison.

Given $X=(Z,U_1,\cdots,U_{k-1})$, $X$ is generated by
$Z\sim unif(-0.5,0.5),$ $U_1=(1+2Z)+e_1,$ $U_2=(1+2Z)+e_2$,
$U_3=(-1+Z)^2+e_3$, and independent $e_j\sim unif(-0.5,0.5)$, for $j=1,2,3$. The other
variables $U_j's$ are generated as: when $3<j<=9$, $U_j=|1+1/(11-j)Z|-|1+1/j\epsilon|$; when $9<j<=19$, $U_j=|1+1/(21-j)Z|-|1+1/j\epsilon|$; and
$U_j=|1+1/(31-j)Z|-|1+1/j\epsilon|$ for $j>19$, $\epsilon \sim unif(-0.5,0.5)$.
We consider the following models in high dimensional setting.
\begin{itemize}
\item { Model 5 (r=1 with uniformed $|Z\cap V^{\top}X|=0$):}
    \begin{eqnarray*}
 Y(1)=\beta_3^{\top}X+\gamma_3Z  U_1  U_2  U_3+\nu,~Y(0)=0~\mbox{and} ~p_5(X)=\Lambda(1+V_3^{\top}X).
\end{eqnarray*}

\item { Model 6 (r=2 with uniformed $|Z\cap V^{\top}X|=1$):}
\begin{eqnarray*}
 Y(1)=\beta_3^{\top}X+\gamma_3Z  U_1  U_2  U_3+\nu,~Y(0)=0,~\mbox{and}~¡¡ p_6(X)=\Lambda(g(\tilde V_3^{\top}X)).
\end{eqnarray*}
\end{itemize}
As for the propensity score, we set  $$			 V_3^{\top}=(\overbrace{{-1,\cdots,-1}}^{5},\overbrace{{0,\cdots,0}}^{5},
\overbrace{{1,\cdots,1}}^{10})/\sqrt{20},\\
\tilde \alpha=(0,\overbrace{{-1,\cdots,-1}}^{4},\overbrace{{0,\cdots,0}}^{5},
\overbrace{{1,\cdots,1}}^{10})/\sqrt{19},$$ and $g(\tilde
V_3^{\top}X)={(1+\tilde \alpha^{\top}X)}/{(1+Z)}$ with   $dim(\tilde V_3^{\top}X)=r=2$, while $|Z\cap V_3^{\top}X|=1$. 
{In high dimensional setting, we only consider the nonlinear model where the parameters are set as $\beta_3^{\top}=(0,\ldots,0)$ and $r_3=1$.
}

The sample size is taken to be $n=500$.
Estimate $\tau(Z)$
 at $Z\in\{-0.4,-0.2,0,0.2,0.4\}$ with 500 simulation realizations.
As for the  bandwidth  choice  we adopt the same rule in
(\ref{sem_band_sele}) of Experiment~1  to have $ h=a n^{ {-1}/(l+4+2r+2\delta_r)}  $ and
  $h_2=a_2 n^{{-1}/(2r+\delta_r)} $. Consider two groups of  $\{a,a_2\},$  $Group~1:$ $\{a=0.55, a_2=0.75\}$ and $Group~2: ~\{a=0.55, a_2=0.69\}$.
For the kernel function in the estimated propensity score, we
also use the Gaussian kernel and higher order kernels derived from it  since the
distribution of $X$ is bounded. All the original simulation results
are reported in
{Table~S.5 in the Supplement} and the relative efficiency results are plotted in Figure~\ref{k=high}.

From the simulation results, we also have the following
 findings.\\{1). The high
dimensionality of $X$ has relatively weak influence on  $IPW$-$S$.  All the
values of $Bias$, $Est\_SD$ and $MSE$  are rather stable and the values of
$P_{\pm1.645}$ are closer to the nominal value $0.05$  as the dimension of
$X$ goes from $4$ to $20$, especially in the case of $|Z \cap V^{\top}X|=1$. } This is very informative because it implies that
$IPW$-$S$ can greatly avoid the curse of dimensionality due to its dimension
reduction structure.\\  2). $IPW$-$S$ not only shows its superiority in dealing
with the curse of dimensionality, but also inherits the efficiency
superiority of $IPW$-$N$ in low-dimensional cases. Under such high dimensional
scenarios, the values of $Est\_SD$ and $MSE$ of $IPW$-$S$ are smaller than
those of the parametric competitors in some cases even $|Z \cap V^{\top}X|=0$. When in Model~6, $IPW$-$S$ is uniformly more efficient than $IPW$-$P$. This is consistent with the theoretical results in Corollary~\ref{cor_sem} as  in the model $|Z\subseteq V^{\top}X|=1$, $IPW$-$S$ is asymptotically more efficient than $IPW$-$P$.  \\

 \begin{flushleft}
 \begin{figure}[H]
 \flushleft
  \includegraphics[scale=0.43]{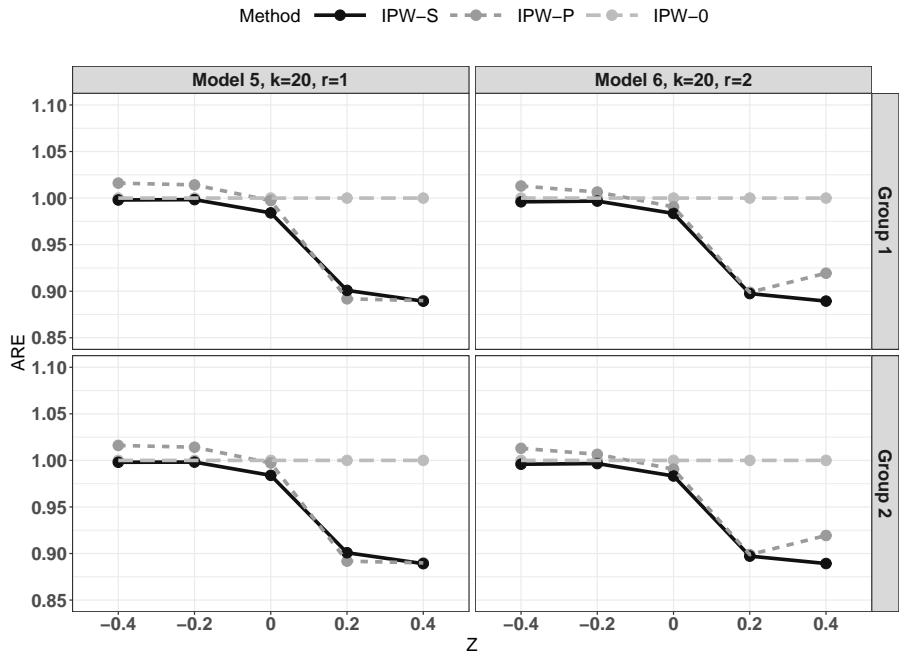}\\
\caption{The asymptotic relative efficiency(ARE) about $Est\_SD$ against that of $IPW$-$O$ under high dimensional setting.
}\label{k=high}
\end{figure}
 \end{flushleft}

{
\section{Data Analysis}

In this section, we consider a dataset collected by \cite{ichino2008}, which can be
obtained from the internet.\footnote{ The data is publicly available at
\url{http://qed.econ.queensu.ca/jae/2008-v23.3/ichino-mealli-nannicini/}}
We
apply the proposed method to estimate the $CATE$ function to
investigate the treatment effect of temporary work assignment ($TWA$) on permanent employment over
worker's age.


First introduce some details and setting about the
dataset. Restricting the sample to Tuscany and aged 17-39,  the resulting
sample size is $n=901$, $281$ of which were on a $TWA$ during the first semester of 2001. That is, the binary treatment variable $D=0,1$ means that the individual was not on or was on a $TWA$ during the first six mouths of 2001. The
outcome $Y$ here is a dummy variable: $Y=1$ if the subject is permanently employed at the end of 2002, and $Y=0$ otherwise. Choose $X_1$ as
  the worker's age and a  set of 25 covariates as $X$ adopted by \cite{ichino2008} to guarantee the unconfoundedness assumption.
   The set of covariates is about demographic characteristics, family background, educational achievements and work experience (See Table~1 in \cite{ichino2008}). This dataset was first analyzed by \cite{ichino2008}, who estimated the parameter $ATT=E(Y_1-Y_0\mid D)$ and showed that $TWA$ can increase the probability of getting a permanent employment. \cite{ichino2008}  pointed out that the $TWA$ effect is heterogeneous for the individuals older than 30 and younger than 30.

    In order to catch more specific heterogeneity of the $TWA$  effect across individuals' age, we  estimate the $CATE$ function $\tau(Z)$ in the interval between ages 20 and 35.
   As the number of covariates is large (= $30$),  we then use a semiparametric single-index
model to estimate the propensity score such that the dimensionality problem and model
misspecification problem can be greatly alleviated. {Given that $D \perp X \mid \alpha^{\top}X$, we can get the $IPW$-$S$ and pointwise confidence band of $\tau(Z)$ by carrying out the estimation procedure proposed in subsection~2.3.}
As for nonparametric estimation part, we use the Gaussian kernel and choose the bandwidths to be $h=0.85\times\hat \sigma_1n^{-1/9}=2.22$, and
$h_2=1.15\times\hat\sigma_dn^{-1/3}=0.04 \ll h,$ where $\hat\sigma_1=\sqrt{var(x_1)}$
and $\hat\sigma_d=\sqrt{var(\alpha^{\top}X)}$.
 We also estimate $IPW$-$P$ as a benchmark to analyse the TWA effect over worker's age.

Figure~\ref{scate_pic} presents the results of  $IPW$-$S$ and $IPW$-$P$ as a function
of worker's age in the range of 20 to 35 years old, which can be regarded as an extension of \cite{ichino2008}  in a certain sense. Furthermore, the 95\% pointwise confidence band of $IPW$-$S$ and $IPW$-$P$  have been also reported in Figure~\ref{scate_pic}.
 There are several points
we want to highlight: 1). {\it both $IPW$-$S$ and $IPW$-$P$ suggest that, from  age $20$ to $35$, a TWA assignment uniformly increases the probability of finding a stable job with the range  roughly between $0.05$ and $0.35$}. It means that if a worker with a TWA experience would more likely to get a permanent job. This finding is in accordance with, but extends the conclusion of \cite{ichino2008}.  2). {\it
The trend of $CATE(x_1)$ varies with worker's age and has two peaks.} From
Figure~\ref{scate_pic}, we can also find that there are two peaks at around age $24$ and age $32$, while the trough appears at around age $29$. That implies the TWA experience has different effect for the workers older than $29$ and under $29$, which was also similarly discussed by \cite{ichino2008}. However, comparing the details in the curves of $IPW$-$S$ and $IPW$-$P$, the effect of TWA on finding a stable job for the subpopulation aged under $29$ is greater than the ones older $29$ in the $IPW$-$S$ curve, while things are opposite in the $IPW$-$P$ curve. It seems that the $IPW$-$S$ curve provides a more reasonable explanation on the  effect of TWA: younger individuals receiving TWA could have better chance to get a stable job than older individuals who need to receive TWA.

\begin{figure}[h]
  \centering
  \includegraphics[scale=0.4]{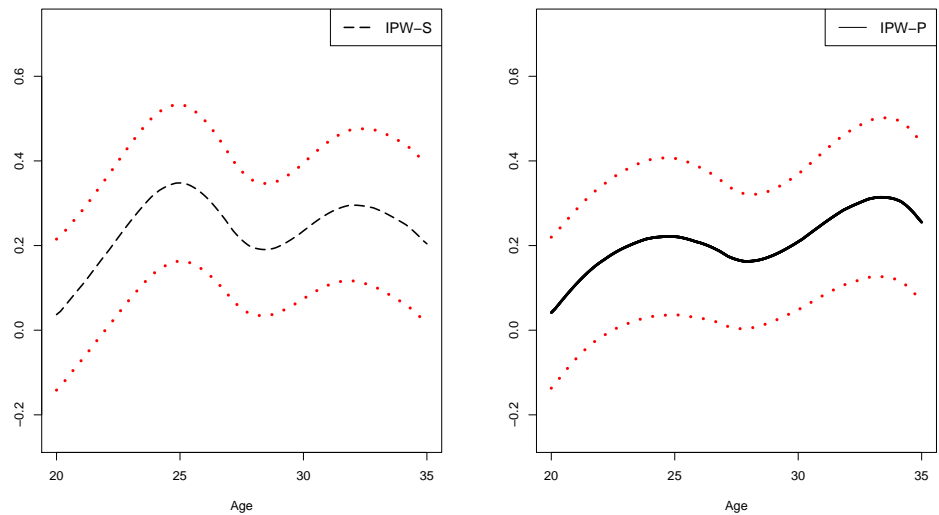}\\
\caption{The curves of conditional average treatment effects (CATE) over worker's age with the 95\% pointwise confidence band.
}\label{scate_pic}
\end{figure}

}

\section {Conclusion}
In this paper, we propose an  estimation ($IPW$-$S$) of conditional
average treatment effect with semiparametric propensity score and investigate its asymptotic properties {which can be used to construct  pointwise confidence intervals}.
We give a relatively complete picture about the asymptotic efficiency of different estimators
 with nonparametric, parametric and true propensity score when model is correctly specified.  Further,
  when
 the dimension of covariates is high, by the numerical studies, we demonstrate the advantages of
 $IPW$-$S$  in alleviating the curse of dimensionality and inheriting the theoretical superiority of $IPW$-$N$ in estimation efficiency.
{But a challenging topic is how to  develop a good uniform confidence band of  the whole function $\tau(z)$ although the Bonferroni confidence band could be applied.} Further, a research topic is about the situation that not all of the covariates are
important for propensity score. Thus, by incorporating variable selection,
we can simultaneously identify important confounders and guarantee the
unconfoundedness assumption. The dimension reduction and variable selection have been investigated by, say, \cite{Ma2019}  for the model under sparsity structure.
This topic is also related to variable selection and thus we will try to
have a computationally inexpensive algorithm for this purpose and study its asymptotic behaviours. Another topic is about the model misspecification even when the semiparametric model is used. We will study the relevant asymptotic behaviours in the near future.

\section*{Acknowledgement}
The authors' s research was supported by grants from NSFC grants (NSFC11671042, NSFC11601227) and the University Grants Council of Hong Kong.

\section*{Supplementary material}
\label{SM}
The supplementary file covers  the detailed proofs to Theorems and Corollaries.

\section*{  Appendix: Technical conditions}

The following regularity conditions are required to get the theoretical results.

\begin{itemize}
\item[(C1)] (Strong ignorability)

 (i) Unconfoundedness: $(Y(0),Y(1))\perp D\mid X$.

 (ii) Common support: For some very small $c>0$, $c<p(X)<1-c$.

\item[(C2)] (on distribution):

(i) The set $\chi$ that is the support of the $k$-dimensional covariate
vector $X$
 is a Cartesian product of compact intervals.\\
(ii) The density function of $Z$, $f_1(Z)$,  and the  density function of
     $X$, are bounded away from zero and infinity and $s \geq r$ times
    continuously differentiable.
\item[(C3)](Conditional moments and smoothness)

(i) $\sup_{x\in\chi}E[Y(j)^2\mid X=x]<\infty$ for $j=0,1;$

 (ii) the functions $m_j(V^{\top}X)=E[Y(j)\mid V^{\top}X ]$, $j=0,1$ are
$s \geq r$ times
    continuously differentiable.

\item[(C4)](on kernel function)

 (i)  $\mathcal{L}(u)$ is a kernel of order $s_1$, is symmetric around zero,
    has finite support $[-1,1]^{\tilde k}$, and is continuously differentiable.

(ii) $\mathcal{H}(u)$ is a kernel of order $s_2$, is symmetric around zero,
    has finite support $[-1,1]^r$, and is continuously differentiable.

 (iii) $K(u)$ is of order $s$, is symmetric around zero,  and
is $s$ times continuously differentiable.
\item[(C5)](on bandwidths)

(i)$h\rightarrow 0$, $nh^l\rightarrow \infty, nh^{2s+l}\rightarrow 0$.

(ii)$h_1, h_2\rightarrow 0$, $log(n)/(nh_2^{r+s_2})\rightarrow 0$, and $log(n)/(nh_1^{\tilde k+s_1})\rightarrow 0.$

(iii)$h_i^{2{s_i}}h^{-2s_{i}-l}\rightarrow 0$, $nh^{l}h_i^{2s_i} \rightarrow 0$, $i=1,2$.
\item[(C6)](on dimension reduction structure)  the dimension of $V$, r, is given and $\hat  V -
    V=O_p(n^{-1/2}).$
\end{itemize}

Recall  the definition of high order kernel in the literature. We say a
function g: $R^r\rightarrow R$ is a kernel of order $s$ if it integrates to
one over $R^r$, and $\int u^{p_1}\cdots u^{p_r}g(u)du=0$ for all nonnegative
integers $p_1,\cdots,p_r$ such that $1\leq \sum_{i}p_i<s, $ and it is
nonzero when $ \sum_{i}p_i=s. $

\section*{References}

\bibliography{jspi_ref}
\bibliographystyle{elsarticle-harv}\biboptions{authoryear}

\begin{table}[H]
  \centering
   \caption{The simulation results under Model~1 with
   $|Z\cap X|=1$ but $|Z\cap V^{\top}X|=0$, $\beta_1^{\top}=(0,0)$ and $\gamma_1=1$  }\label{k=2_r=1_G1_s1}
  \renewcommand\arraystretch{0.65}
  \scalebox{0.65}{
\begin{tabular}{ccrrrrrrrrrrr}
    \hline
    \multicolumn{13}{c}{k=2, r=1 under Group 1:  $\{a=0.55,a_1=1.05,a_2=0.75\}$} \\
    \cmidrule{1-7}\cmidrule{8-13}    n     & $Z$   & \multicolumn{1}{c}{-0.4} & \multicolumn{1}{c}{-0.2} & \multicolumn{1}{c}{0} & \multicolumn{1}{c}{0.2} & \multicolumn{1}{c}{0.4} &       & \multicolumn{1}{c}{-0.4} & \multicolumn{1}{c}{-0.2} & \multicolumn{1}{c}{0} & \multicolumn{1}{c}{0.2} & \multicolumn{1}{c}{0.4} \\
           \cmidrule{1-7}\cmidrule{9-13}          &       & \multicolumn{5}{c}{$IPW$-$O$}             &       & \multicolumn{5}{c}{$IPW$-$P$} \\
\cmidrule{2-7}\cmidrule{9-13}    \multirow{11}[6]{*}{500} & Bias  & -0.0262  & 0.0607  & 0.0553  & -0.0278  & -0.1072  &       & -0.0262  & 0.0607  & 0.0553  & -0.0276  & -0.1068  \\
          & $Est_SD$ & 0.0269  & 0.0199  & 0.0190  & 0.0215  & 0.0271  &       & 0.0265  & 0.0197  & 0.0190  & 0.0204  & 0.0238  \\
          & MSE   & 0.0014  & 0.0041  & 0.0034  & 0.0012  & 0.0122  &       & 0.0014  & 0.0041  & 0.0034  & 0.0012  & 0.0120  \\
          & $P_{-1.645}$ & 0.0540  & 0.0440  & 0.0580  & 0.0560  & 0.0560  &       & 0.0520  & 0.0400  & 0.0540  & 0.0540  & 0.0380  \\
          & $P_{1.645}$ & 0.0480  & 0.0460  & 0.0400  & 0.0360  & 0.0540  &       & 0.0520  & 0.0500  & 0.0420  & 0.0400  & 0.0600  \\
\hline
    \end{tabular}%
}%
\end{table}%

\begin{table}[H]
  \centering
   \caption{The simulation results under Model~1 with
   $|Z\cap U|=0$ and $|Z\cap
    V^{\top}X|=0$, $\beta_1^{\top}=(1/2,-1/5)$ and $\gamma_1=0$ }\label{k=2_r=1_G1_s2}
  \renewcommand\arraystretch{0.65}
  \scalebox{0.7}{
    \begin{tabular}{ccrrrrrrrrrrr}
    \hline
    \multicolumn{13}{c}{k=2, r=1 under Group 1:  $\{a=0.55,a_1=1.05,a_2=0.75\}$} \bigstrut\\
   \cmidrule{1-7}\cmidrule{8-13}    n     & $Z$   & \multicolumn{1}{c}{-0.4} & \multicolumn{1}{c}{-0.2} & \multicolumn{1}{c}{0} & \multicolumn{1}{c}{0.2} & \multicolumn{1}{c}{0.4} &       & \multicolumn{1}{c}{-0.4} & \multicolumn{1}{c}{-0.2} & \multicolumn{1}{c}{0} & \multicolumn{1}{c}{0.2} & \multicolumn{1}{c}{0.4} \\
           \cmidrule{1-7}\cmidrule{9-13}          &       & \multicolumn{5}{c}{$IPW$-$O$}             &       & \multicolumn{5}{c}{$IPW$-$P$} \\
\cmidrule{2-7}\cmidrule{9-13}
     \multirow{11}[6]{*}{500} & Bias  & 0.0027  & -0.0006  & 0.0174  & 0.0147  & -0.0504  &       & 0.0021  & -0.0010  & 0.0173  & 0.0148  & -0.0503  \bigstrut[t]\\
          & $Est_SD$ & 0.0297  & 0.0215  & 0.0185  & 0.0187  & 0.0211  &       & 0.0254  & 0.0190  & 0.0174  & 0.0181  & 0.0210  \\
          & MSE   & 0.0009  & 0.0005  & 0.0006  & 0.0006  & 0.0030  &       & 0.0006  & 0.0004  & 0.0006  & 0.0005  & 0.0030  \\
          & $P_{-1.645}$ & 0.0540  & 0.0500  & 0.0420  & 0.0500  & 0.0420  &       & 0.0460  & 0.0520  & 0.0420  & 0.0540  & 0.0440  \\
          & $P_{1.645}$ & 0.0440  & 0.0440  & 0.0440  & 0.0480  & 0.0560  &       & 0.0480  & 0.0400  & 0.0480  & 0.0480  & 0.0560  \bigstrut[b]\\
\hline
    \end{tabular}}
\end{table}%

\begin{table}[H]
  \centering
   \caption{The simulation results under Model~2 with
   $|Z\cap U|=0$ and $|Z\cap
    V^{\top}X|=0$, $\beta_1^{\top}=(0,0)$ and $\gamma_1=1$ }\label{k=1_r=1_G1_s1}
  \renewcommand\arraystretch{0.65}
  \scalebox{0.7}{
    \begin{tabular}{ccrrrrrrrrrrr}
    \hline
    \multicolumn{13}{c}{k=2, r=1 under Group 1:  $\{a=0.55,a_1=1.05,a_2=0.75\}$} \\
  \cmidrule{1-7}\cmidrule{8-13}    n     & $Z$   & \multicolumn{1}{c}{-0.4} & \multicolumn{1}{c}{-0.2} & \multicolumn{1}{c}{0} & \multicolumn{1}{c}{0.2} & \multicolumn{1}{c}{0.4} &       & \multicolumn{1}{c}{-0.4} & \multicolumn{1}{c}{-0.2} & \multicolumn{1}{c}{0} & \multicolumn{1}{c}{0.2} & \multicolumn{1}{c}{0.4} \\
           \cmidrule{1-7}\cmidrule{9-13}          &       & \multicolumn{5}{c}{$IPW$-$O$}             &       & \multicolumn{5}{c}{$IPW$-$P$} \\
\cmidrule{2-7}\cmidrule{9-13}  \multirow{11}[6]{*}{500} & Bias  & -0.0269  & 0.0603  & 0.0550  & -0.0282  & -0.1080  &       & -0.0269  & 0.0603  & 0.0551  & -0.0278  & -0.1072  \\
          & $Est_SD$ & 0.0255  & 0.0192  & 0.0186  & 0.0208  & 0.0262  &       & 0.0252  & 0.0190  & 0.0186  & 0.0199  & 0.0233  \\
          & MSE   & 0.0014  & 0.0040  & 0.0034  & 0.0012  & 0.0123  &       & 0.0014  & 0.0040  & 0.0034  & 0.0012  & 0.0120  \\
          & $P_{-1.645}$ & 0.0520  & 0.0480  & 0.0520  & 0.0540  & 0.0440  &       & 0.0520  & 0.0480  & 0.0540  & 0.0600  & 0.0400  \\
          & $P_{1.645}$ & 0.0520  & 0.0480  & 0.0420  & 0.0420  & 0.0540  &       & 0.0440  & 0.0520  & 0.0480  & 0.0420  & 0.0620  \\
\hline
    \end{tabular}%
    }%
\end{table}%

\begin{table}[H]
  \centering
   \caption{The simulation results under Model~2 with
   $|Z\cap U|=0$ and $|Z\cap
    V^{\top}X|=0$, $\beta_1^{\top}=(1/2,-1/5)$ and $\gamma_1=0$ }\label{k=1_r=1_G1_s2}
  \renewcommand\arraystretch{0.65}
  \scalebox{0.7}{
    \begin{tabular}{ccrrrrrrrrrrr}
    \hline
    \multicolumn{13}{c}{Model~2 k=2, r=1 under Group 1:  $\{a=0.55,a_1=1.05,a_2=0.75\}$} \bigstrut\\
    \cmidrule{1-7}\cmidrule{8-13}    n     & $Z$   & \multicolumn{1}{c}{-0.4} & \multicolumn{1}{c}{-0.2} & \multicolumn{1}{c}{0} & \multicolumn{1}{c}{0.2} & \multicolumn{1}{c}{0.4} &       & \multicolumn{1}{c}{-0.4} & \multicolumn{1}{c}{-0.2} & \multicolumn{1}{c}{0} & \multicolumn{1}{c}{0.2} & \multicolumn{1}{c}{0.4} \\
           \cmidrule{1-7}\cmidrule{9-13}          &       & \multicolumn{5}{c}{$IPW$-$O$}             &       & \multicolumn{5}{c}{$IPW$-$P$} \\
\cmidrule{2-7}\cmidrule{9-13}    \multirow{11}[6]{*}{500} & Bias  & 0.0015  & -0.0014  & 0.0172  & 0.0148  & -0.0504  &       & 0.0015  & -0.0014  & 0.0171  & 0.0147  & -0.0504  \bigstrut[t]\\
          & $Est_SD$ & 0.0259  & 0.0196  & 0.0175  & 0.0182  & 0.0208  &       & 0.0232  & 0.0181  & 0.0168  & 0.0179  & 0.0208  \\
          & MSE   & 0.0007  & 0.0004  & 0.0006  & 0.0005  & 0.0030  &       & 0.0005  & 0.0003  & 0.0006  & 0.0005  & 0.0030  \\
          & $P_{-1.645}$ & 0.0580  & 0.0460  & 0.0520  & 0.0540  & 0.0460  &       & 0.0480  & 0.0480  & 0.0520  & 0.0540  & 0.0440  \\
          & $P_{1.645}$ & 0.0480  & 0.0460  & 0.0540  & 0.0440  & 0.0480  &       & 0.0460  & 0.0580  & 0.0460  & 0.0500  & 0.0460  \bigstrut[b]\\
\hline
    \end{tabular}}
\end{table}%

\begin{table}[H]
  \centering
  \caption{The simulation results under Model~3 with
   $|Z\cap X|=1$ but $|Z\cap
    V^{\top}X|=0$, $\beta_2^{\top}=(1/10,1/\sqrt{2},-1/\sqrt{2},-1/10)$ and $\gamma_2=0$ }\label{k=4_r=1_s1}
    \renewcommand\arraystretch{0.65}
  \scalebox{0.7}{
         \begin{tabular}{ccrrrrrrrrrrr}
    \hline
    \multicolumn{13}{c}{k=4, r=1 under Group 1:  $\{a=0.55,a_1=1.05,a_2=0.75\}$} \\
    \cmidrule{1-7}\cmidrule{8-13}    n     & $Z$   & \multicolumn{1}{c}{-0.4} & \multicolumn{1}{c}{-0.2} & \multicolumn{1}{c}{0} & \multicolumn{1}{c}{0.2} & \multicolumn{1}{c}{0.4} &       & \multicolumn{1}{c}{-0.4} & \multicolumn{1}{c}{-0.2} & \multicolumn{1}{c}{0} & \multicolumn{1}{c}{0.2} & \multicolumn{1}{c}{0.4} \\
           \cmidrule{1-7}\cmidrule{9-13}          &       & \multicolumn{5}{c}{$IPW$-$O$}             &       & \multicolumn{5}{c}{$IPW$-$P$} \\
\cmidrule{2-7}\cmidrule{9-13}     \multirow{11}[6]{*}{500} & Bias  & -0.0271  & 0.0610  & 0.0556  & -0.0280  & -0.1086  &       & -0.0271  & 0.0611  & 0.0557  & -0.0279  & -0.1083  \\
          & $Est_SD$ & 0.0225  & 0.0179  & 0.0176  & 0.0232  & 0.0334  &       & 0.0225  & 0.0179  & 0.0177  & 0.0231  & 0.0329  \\
          & MSE   & 0.0012  & 0.0040  & 0.0034  & 0.0013  & 0.0129  &       & 0.0012  & 0.0040  & 0.0034  & 0.0013  & 0.0128  \\
          & $P_{-1.645}$ & 0.0360  & 0.0380  & 0.0540  & 0.0420  & 0.0400  &       & 0.0340  & 0.0420  & 0.0520  & 0.0480  & 0.0460  \\
          & $P_{1.645}$ & 0.0480  & 0.0420  & 0.0460  & 0.0500  & 0.0440  &       & 0.0460  & 0.0440  & 0.0420  & 0.0560  & 0.0420  \\
\hline
    \end{tabular}%
}
\end{table}%

\begin{table}[H]
  \centering
  \caption{The simulation results under Model~3 with
   $|Z\cap X|=1$ but $|Z\cap
    V^{\top}X|=0$, $\beta_2^{\top}=(1/10,1/\sqrt{2},-1/\sqrt{2},-1/10)$ and $\gamma_2=0$ }\label{k=4_r=1_s2}
    \renewcommand\arraystretch{0.65}
  \scalebox{0.7}{
    \begin{tabular}{ccrrrrrrrrrrr}
    \hline
    \multicolumn{13}{c}{k=4, r=1 under Group 1:  $\{a=0.55,a_1=1.05,a_2=0.75\}$} \bigstrut\\
     \cmidrule{1-7}\cmidrule{8-13}    n     & $Z$   & \multicolumn{1}{c}{-0.4} & \multicolumn{1}{c}{-0.2} & \multicolumn{1}{c}{0} & \multicolumn{1}{c}{0.2} & \multicolumn{1}{c}{0.4} &       & \multicolumn{1}{c}{-0.4} & \multicolumn{1}{c}{-0.2} & \multicolumn{1}{c}{0} & \multicolumn{1}{c}{0.2} & \multicolumn{1}{c}{0.4} \\
           \cmidrule{1-7}\cmidrule{9-13}          &       & \multicolumn{5}{c}{$IPW$-$O$}             &       & \multicolumn{5}{c}{$IPW$-$P$} \\
\cmidrule{2-7}\cmidrule{9-13}     \multirow{11}[6]{*}{500} & Bias  & 0.0333  & 0.0042  & -0.0019  & -0.0067  & -0.0228  &       & 0.0332  & 0.0043  & -0.0018  & -0.0066  & -0.0228  \bigstrut[t]\\
          & $Est_SD$ & 0.0317  & 0.0264  & 0.0236  & 0.0236  & 0.0279  &       & 0.0300  & 0.0247  & 0.0222  & 0.0227  & 0.0273  \\
          & MSE   & 0.0021  & 0.0007  & 0.0006  & 0.0006  & 0.0013  &       & 0.0020  & 0.0006  & 0.0005  & 0.0006  & 0.0013  \\
          & $P_{-1.645}$ & 0.0440  & 0.0520  & 0.0560  & 0.0460  & 0.0500  &       & 0.0480  & 0.0540  & 0.0500  & 0.0560  & 0.0540  \\
          & $P_{1.645}$ & 0.0480  & 0.0440  & 0.0400  & 0.0580  & 0.0480  &       & 0.0440  & 0.0440  & 0.0380  & 0.0500  & 0.0560  \bigstrut[b]\\
\hline
    \end{tabular}}
\end{table}%

\begin{table}[H]
  \centering
   \caption{The simulation results under Model~4 with k=4, r=2,
   $|Z\cap X|=1$, $|Z\cap
    V^{\top}X|=1$, $\beta_2^{\top}=(0,0,0,0)$ and $\gamma_2=1$ }\label{k=4_r=2_s1}
  \renewcommand\arraystretch{0.65}
  \scalebox{0.7}{
    \begin{tabular}{ccrrrrrrrrrrr}
    \hline
    \multicolumn{13}{c}{k=4, r=2 under Group 1: $\{a=0.55,a_1=1.05,a_2=0.75\}$} \\
   \cmidrule{1-7}\cmidrule{8-13}    n     & $Z$   & \multicolumn{1}{c}{-0.4} & \multicolumn{1}{c}{-0.2} & \multicolumn{1}{c}{0} & \multicolumn{1}{c}{0.2} & \multicolumn{1}{c}{0.4} &       & \multicolumn{1}{c}{-0.4} & \multicolumn{1}{c}{-0.2} & \multicolumn{1}{c}{0} & \multicolumn{1}{c}{0.2} & \multicolumn{1}{c}{0.4} \\
           \cmidrule{1-7}\cmidrule{9-13}          &       & \multicolumn{5}{c}{$IPW$-$O$}             &       & \multicolumn{5}{c}{$IPW$-$P$} \\
\cmidrule{2-7}\cmidrule{9-13}     \multirow{11}[6]{*}{500} & Bias  & -0.0260  & 0.0621  & 0.0565  & -0.0277  & -0.1087  &       & -0.0258  & 0.0623  & 0.0567  & -0.0274  & -0.1082  \\
          & $Est_SD$ & 0.0247  & 0.0201  & 0.0208  & 0.0304  & 0.0469  &       & 0.0247  & 0.0201  & 0.0207  & 0.0279  & 0.0416  \\
          & MSE   & 0.0013  & 0.0043  & 0.0036  & 0.0017  & 0.0140  &       & 0.0013  & 0.0043  & 0.0036  & 0.0015  & 0.0134  \\
          & $P_{-1.645}$ & 0.0460  & 0.0500  & 0.0500  & 0.0440  & 0.0400  &       & 0.0380  & 0.0520  & 0.0540  & 0.0480  & 0.0460  \\
          & $P_{1.645}$ & 0.0500  & 0.0540  & 0.0500  & 0.0500  & 0.0460  &       & 0.0460  & 0.0400  & 0.0500  & 0.0560  & 0.0500  \\
\hline
    \end{tabular}%
}
\end{table}%

\begin{table}[H]
  \centering
   \caption{The simulation results under Model~4 with k=4, r=2,
   $|Z\cap X|=1$ and $|Z\cap
    V^{\top}X|=1$, $\beta_2^{\top}=(1/10,1/\sqrt{2},-1/\sqrt{2},-1/10)$ and $\gamma_2=0$ }\label{k=4_r=2_s2}
  \renewcommand\arraystretch{0.65}
  \scalebox{0.7}{
    \begin{tabular}{ccrrrrrrrrrrr}
    \hline
    \multicolumn{13}{c}{k=4, r=2 under Group 1: $\{a=0.55,a_1=1.05,a_2=0.75\}$} \bigstrut\\
    \cmidrule{1-7}\cmidrule{8-13}    n     & $Z$   & \multicolumn{1}{c}{-0.4} & \multicolumn{1}{c}{-0.2} & \multicolumn{1}{c}{0} & \multicolumn{1}{c}{0.2} & \multicolumn{1}{c}{0.4} &       & \multicolumn{1}{c}{-0.4} & \multicolumn{1}{c}{-0.2} & \multicolumn{1}{c}{0} & \multicolumn{1}{c}{0.2} & \multicolumn{1}{c}{0.4} \\
           \cmidrule{1-7}\cmidrule{9-13}          &       & \multicolumn{5}{c}{$IPW$-$O$}             &       & \multicolumn{5}{c}{$IPW$-$P$} \\
\cmidrule{2-7}\cmidrule{9-13}    \multirow{11}[6]{*}{500} & Bias  & 0.0339  & 0.0050  & -0.0009  & -0.0057  & -0.0224  &       & 0.0343  & 0.0053  & -0.0008  & -0.0055  & -0.0222  \bigstrut[t]\\
          & $Est_SD$ & 0.0364  & 0.0298  & 0.0278  & 0.0296  & 0.0347  &       & 0.0356  & 0.0295  & 0.0277  & 0.0295  & 0.0347  \\
          & MSE   & 0.0025  & 0.0009  & 0.0008  & 0.0009  & 0.0017  &       & 0.0024  & 0.0009  & 0.0008  & 0.0009  & 0.0017  \\
          & $P_{-1.645}$ & 0.0520  & 0.0480  & 0.0560  & 0.0460  & 0.0480  &       & 0.0520  & 0.0460  & 0.0580  & 0.0380  & 0.0560  \\
          & $P_{1.645}$ & 0.0440  & 0.0520  & 0.0420  & 0.0520  & 0.0440  &       & 0.0380  & 0.0440  & 0.0400  & 0.0540  & 0.0440  \bigstrut[b]\\
\hline
    \end{tabular}}

\end{table}%

\begin{landscape}
\begin{table}[H]
 \flushleft
   \caption{The simulation results under high dimensional setting}\label{k=high_unif}
  \renewcommand\arraystretch{0.7}
  \scalebox{0.75}{
       \begin{tabular}{ccrrrrrcrrrrrcrrrrr}
    \hline
    \multicolumn{19}{c}{ dim(X)=20 with $|Z\cap V^{\top}X|=0$: $r=1$} \\
    \hline
          &       & \multicolumn{5}{c}{$IPW$-$O$}             &       & \multicolumn{5}{c}{$IPW$-$P$}             &       & \multicolumn{5}{c}{$IPW$-$S$} \\
\cmidrule{3-7}\cmidrule{9-13}\cmidrule{15-19}    Group & $Z$   & \multicolumn{1}{c}{-0.4} & \multicolumn{1}{c}{-0.2} & \multicolumn{1}{c}{0} & \multicolumn{1}{c}{0.2} & \multicolumn{1}{c}{0.4} &       & \multicolumn{1}{c}{-0.4} & \multicolumn{1}{c}{-0.2} & \multicolumn{1}{c}{0} & \multicolumn{1}{c}{0.2} & \multicolumn{1}{c}{0.4} &       & \multicolumn{1}{c}{-0.4} & \multicolumn{1}{c}{-0.2} & \multicolumn{1}{c}{0} & \multicolumn{1}{c}{0.2} & \multicolumn{1}{c}{0.4} \\
\cmidrule{1-7}\cmidrule{9-13}\cmidrule{15-19}    \multirow{5}[2]{*}{Group~1} & Bias  & -0.0275  & 0.0598  & 0.0554  & -0.0271  & -0.1064  &       & -0.0273  & 0.0601  & 0.0556  & -0.0270  & -0.1062  &       & -0.0254  & 0.0612  & 0.0538  & -0.0349  & -0.1209  \\
          & $Est\_SD$ & 0.0244  & 0.0201  & 0.0211  & 0.0296  & 0.0474  &       & 0.0248  & 0.0204  & 0.0210  & 0.0264  & 0.0422  &       & 0.0243  & 0.0201  & 0.0208  & 0.0267  & 0.0422  \\
          & MSE   & 0.0014  & 0.0040  & 0.0035  & 0.0016  & 0.0136  &       & 0.0014  & 0.0040  & 0.0035  & 0.0014  & 0.0131  &       & 0.0012  & 0.0041  & 0.0033  & 0.0019  & 0.0164  \\
          & $P_{-1.645}$ & 0.0420  & 0.0520  & 0.0420  & 0.0380  & 0.0420  &       & 0.0580  & 0.0500  & 0.0420  & 0.0460  & 0.0540  &       & 0.0520  & 0.0560  & 0.0420  & 0.0440  & 0.0500  \\
          & $P_{1.645}$ & 0.0480  & 0.0560  & 0.0540  & 0.0600  & 0.0540  &       & 0.0440  & 0.0540  & 0.0500  & 0.0540  & 0.0540  &       & 0.0460  & 0.0540  & 0.0620  & 0.0520  & 0.0560  \\
\hline
    \end{tabular}%
}
\end{table}%
\end{landscape}

\end{document}